\newtheorem{theo}{Theorem}
\newtheorem{prop}[theo]{Proposition}
\newtheorem{Lemme}[theo]{Lemma}
\newtheorem{cor}[theo]{Corollary}

\newcommand {\pare}[1] {\left( {#1} \right)}
\newcommand {\cro}[1] {\left[ {#1} \right]}
\newcommand {\bra}[1] {\left\langle {#1} \right\rangle}
\newcommand {\acc}[1] {\left\{ {#1} \right\}}
\newcommand {\nor}[1] { \left\| {#1} \right\|}
\newcommand {\va}[1] {\left| {#1} \right|}

\newcommand {\sous}[2] {\begin{array}[t]{c}
                          #1 \\[-2mm]
              \scriptstyle{#2}
            \end{array}}
\newcommand {\sur}[2] { \stackrel {\scriptstyle{#1}}{#2}}
\newcommand {\inte}[1] {\begin{array}[b]{c}
			 \hspace{4pt} \scriptstyle{\circ} \\[-3mm]
				{#1}
            		\end{array}}
\newcommand {\refeq}[1] {(\ref{#1})}
\newcommand {\Frac}[2] {\displaystyle{\frac {#1}{#2} }}
\newcommand {\limisup}[1] {\displaystyle{\limsup_{#1}}}
\newcommand {\limiinf}[1] {\displaystyle{\liminf_{#1}}}
\newcommand {\limi}[1] {\displaystyle{\lim_{#1}}}
\newcommand {\som}[2] {\displaystyle{\sum_{#1}^{#2}}}
\newcommand {\pro}[2] {\displaystyle{\prod_{#1}^{#2}}}

\newcommand{\Inf}{\mathop{\mathrm{Inf}}}
\newcommand{\Min}{\mathop{\mathrm{Min}}}
\newcommand{\Max}{\mathop{\mathrm{Max}}}
\newcommand{\Sup}{\mathop{\mathrm{Sup}}}

\def \R {\mathbb{R}} 
\def \P  {\mathbb{P}} 
\def \P  {\mathbb{P}} 

\def \Z  {\mathbb{Z}}
\def \Q {\mathbb{Q}}
\def \N {\mathbb{N}}
\def \D {\mathbb{D}}
\def \ind {\hbox{ 1\hskip -3pt I}}

\def \LL {{\mathcal L}}
\def \BB {{\mathcal B}}
\def \DD {{\mathcal D}}
\def \CC {{\mathcal C}}
\def \II {{\mathcal I}}
\def \JJ {{\mathcal J}}
\def \XX {{\mathcal X}}
\def \GG {{\mathcal G}} 
\def \MM {{\mathcal M}}
\def \AA {{\mathcal A}}
\def \NN {{\mathcal N}}
\def \UU {{\mathcal U}}

\documentclass[12pt]{article}
\usepackage{amssymb}
\usepackage{eucal}
\begin{document}

\vspace*{.2cm}
\begin{center}
{\Large{\bf
 QUENCHED LARGE DEVIATIONS  FOR DIFFUSIONS IN A RANDOM GAUSSIAN
SHEAR FLOW DRIFT.
}}
\vspace{5mm}\\
 {\large{\bf Amine Asselah}} \footnote{E-mail: asselah@cmi.univ-mrs.fr},
 {\large{\bf Fabienne Castell}} \footnote{E-mail: castell@cmi.univ-mrs.fr}.\\
 \vspace{.5mm}
 Laboratoire d'Analyse, Topologie
 et Probabilit\'es. CNRS UMR 6632. \\
 CMI. Universit\'e de Provence. \\
 39 rue Joliot Curie.  \\
 13453 Marseille Cedex 13. FRANCE.
\end{center}
\vspace{1cm}

{\bf Abstract.}
{\small We prove a full large deviations principle  in large time, for a
diffusion process with random drift
$X_t = W_t + \int_{0}^{t} V \pare{{X_s}} \, ds$,
where $V$ is a centered Gaussian shear flow random field independent of the 
Brownian $W$. The large deviations principle is established in a ``quenched'' 
setting, i.e.  is  valid almost surely in the randomness  of $V$.
}
\vspace{.5cm}


Large deviations. Parabolic Anderson model.

\vspace{1cm}
\section{Introduction.}
\indent

In this paper, we investigate  large deviations properties for
diffusions $(X_t, t \geq 0)$ with random drift, solving 
\begin{equation} 
\label{sde}
X_t = W_t + \int_{0}^{t} V \pare{X_s} 
\, ds \,\, ,
\end{equation} 
where $W$ is a standard Brownian motion in $\R^2$, and
$V$ is a centered stationary solenoidal (i.e. such that
$\mbox{div}(V)=0$) Gaussian field on $\R^2$,
independent of $W$. 

Such a process is a model for diffusion in
an incompressible  turbulent flow. As such, it has been discussed 
thoroughly both in the physics and ma\-the\-matics 
literature (see for instance 
\cite{avellaneda-majda90,avellaneda-majda92,carmona,carmona-xu,landim-olla-yau,olla}).  These papers deal with the long
time behavior of the process $X$. More precisely they
investigate the link between the properties of the 
random drift $V$, and the convergence in law of $X_t$
when $t$ goes to infinity.

The model we are working on in this paper, is a very particular case
of \refeq{sde}, since $V$ is assumed to be a {\bf shear flow}, i.e.
\begin{equation} 
\label{drift} 
\forall x \in \R^2 \, ,  \,\,\,  x = (x_1, x_2) \, ,  \,\,\,\, 
V(x) = (0, v(x_1) )\,\, .
\end{equation} 
$(v(x_1), x_1 \in \R) $ is a centered Gaussian 
field, with covariance
$K(x_1 - x_1') \triangleq 
\bra{v(x_1) v(x_1')}$. This model has the advantage 
of being easy to handle, since in the shear flow situation, the two
coordinates $(X_{1,T}, X_{2,T})$ of $X_T$ are just
\begin{equation}
\label{coord}
\left\{ \begin{array}{ll}
X_{1,T} & = W_{1,T} \\
X_{2,T} & = W_{2,T} + \int_0^T v(W_{1,s}) \, ds
\end{array}
\right.
\end{equation} 

>From the viewpoint of
the central limit theorem, this model has been studied in 
\cite{avellaneda-majda90,avellaneda-majda92}, where it is proved that when 
the covariance function $K$ decays  sufficiently slowly at infinity,
the second coordinate $X_{2,T}$ of $X$ exhibits a super-diffusive
behavior, i.e. for some parameter $\alpha  > 1/2$  (related to the decay
of correlation), 
$ \frac{1}{T^{\alpha }} X_{2,T}$ converges in law  when $T \rightarrow \infty$.

In \cite{castell-pradeilles}, the annealed large deviations of the
Gaussian shear flow model \refeq{sde} \refeq{drift}
are established. The result is the 
following. Let $\P$ denote the annealed law, that is the law of $X$
integrated over the randomnesses of $V$ and $W$. Then, for all Borel set $A$ of
$\R$, with closure $\bar{A}$, and interior $\sur{\, \, \circ}{A}$.
\begin{equation} 
\begin{array}{ll}
- \begin{array}[t]{c} 
  \Inf \\[-10pt]
   \scriptstyle{x}  \in \hspace{-2mm} \inte{\scriptstyle{A}}
  \end{array}  L(x) &
\leq  \liminf_{T \rightarrow \infty} \frac{1}{T} 
\log \P\cro{ \frac{1}{T^{3/2}} X_{2,T} \in A} \\
& \leq \limsup_{T \rightarrow \infty} \frac{1}{T} 
\log \P\cro{ \frac{1}{T^{3/2}} X_{2,T} \in A}
\leq - \sous{\Inf}{x \in \bar{A}} L(x) \, . 
\end{array} 
\end{equation}
The rate function $L$ is continuous, with compact level sets
and  has a unique zero  at the origin. Note that the super-diffusive
scaling $T^{3/2}$ does not depend on the decay of correlation,
but is intimately linked with  the choice of Gaussian statistics 
for the drift $V$.

  We study here the large deviations of the Gaussian shear flow
model in a quenched setting, i.e. almost surely in the environment
$V$. Our main result states that there exists 
a convex  deterministic rate function  $\JJ$ 
such that a.s. in $V$ and for all Borel set $A$ of $\R$, 
\begin{equation}
\label{resultat}
\begin{array}{ll} 
- \begin{array}[t]{c} 
  \Inf \\[-10pt]
   \scriptstyle{x}  \in \hspace{-2mm} \inte{\scriptstyle{A}}
  \end{array}  \JJ(x) &
\leq  \liminf_{T \rightarrow \infty} \frac{1}{T}
\log \P\cro{ \frac{1}{T \sqrt{\log(T)}} X_{2,T} \in A} \\
 & \leq  \limsup_{T \rightarrow \infty} \frac{1}{T} 
\log \P\cro{ \frac{1}{T \sqrt{\log(T)}} X_{2,T} \in A}
\leq - \sous{\Inf}{x \in \bar{A}} \JJ(x) \, . 
\end{array} 
\end{equation}   
Note that in the scaling $T \sqrt{\log(T)}$, the Brownian part
$W_2$ does not play any role in the large deviations result, 
and \refeq{resultat} states  a  large deviations principle  for 
$Y_T \triangleq  \frac{1}{T \sqrt{\log(T)}} \int_0^T v(B_s) \,ds$, where
$B$ is a Brownian motion independent of $v$.
Here again, the super-diffusive scaling does not depend
on the decay  of correlation, but on the choice of
the Gaussian law  for $v$. This scaling is related to
the order of magnitude of a Gaussian field on a box of size
$T$.  Indeed, with probability of order 
$\exp(-RT)$ ($R$ large), the Brownian motion $(B_s, 0 \leq s\leq T)$
stays in a ball of radius $\sqrt{R} T$, so that in the study of 
the large deviations of $Y_T$, we can restrict ourselves to 
trajectories confined to such balls with $R$ large enough.
 So the effect of the scaling
$\sqrt{\log(T)}$ is to deal with a bounded integrand $v/\sqrt{\log(T)}$. 

The  large deviations upper bound is obtained using the G\"artner-Ellis
method, i.e. by considering the quenched behavior of the Laplace transform
\begin{equation} 
\label{laplace}
\begin{array}{ll}
\Lambda_T(\alpha) & = E_0 \cro{ \exp(\alpha T Y_T); \tau_{RT} \geq T}
\\[2mm]
& = E_0 \cro{ \exp \pare{ \alpha \int_0^T \frac{v(B_s)}{\sqrt{\log(T)}} \, ds} 
; \tau_{RT} \geq T} \,\, .
\end{array} 
\end{equation} 
In expression \refeq{laplace},
 $E_0$ denotes the expectation with respect to $B$, assuming that
$B_0=0$, and $\tau_{RT}$ is the first exit time of $B$ from
the interval $I_{RT}=]-RT;RT[$.
As usual, we are led to look at the a.s limit when $T \rightarrow \infty$,
of the principal eigenvalue of the random operator
$\LL(f) = -\frac{1}{2} f'' - \alpha \frac{v}{\sqrt{\log(T)}} f$,
with Dirichlet conditions on the boundary of $I_{RT}$,
\begin{equation}	
\label{vpp}
\begin{array}{l}
\lambda \pare{\alpha v/\sqrt{\log(T)}, B(0,RT)}  
\\ 
 = \Inf  \left\{ \frac{1}{2} \int (f')^2(x) \, dx
	- \int \frac{\alpha v(x)}{\sqrt{\log(T)}} 
		f^2(x) dx : f \in \CC^{\infty}_c(B(0,RT)), \right.
\\ \hspace{8cm}	 \left.
		\int f^2(x) dx =1 
\right\} \, .
\end{array} 
\end{equation}  
Following the image popularized by A.S Sznitman in \cite{sznitman},
the main contribution comes from ``the regions where the 
eigenvalue is small''. 
Thus, a key argument in the study of 
$\lambda \pare{\alpha v/\sqrt{\log(T)}, B(0,RT)}$
is a lemma borrowed from \cite{gartner-konig}, 
which asserts that this principal eigenvalue
is   comparable with $\Min_{i} \lambda(\alpha v/\sqrt{\log(T)},Q_i)$, where
$Q_i$ are balls of fixed size covering  $I_{RT}$. This comparison
enables one to show that $v$-a.s, $\Lambda(\alpha) 
\triangleq \lim_{T \rightarrow \infty} 
(\log \Lambda_T(\alpha))/T $ exists, and is deterministic. The upper
bound is thus obtained with a rate functional $\JJ$ 
which is the Legendre transform of
$\Lambda$.

  On the opposite direction, a
 first lower bound is obtained using a specific strategy 
for the path of the Brownian motion: we force it  to go
``fast'' to a region where the field $v$ has a ``high'' peak, and to
remain there until time $T$. The  rate function $\II_1$ obtained 
in this way, has a Legendre transform which coincides
with $\JJ$. Thus, if $\II_1$ were convex, then $\II_1=\JJ$.
However, we could not prove convexity of $\II_1$.
We overcome this problem by adopting the following strategy.
We imagine a sequence of scenarii: the $n$-th one corresponds to
partitioning $[0,T]$ into  $n$  time intervals, in each of which  
the Brownian motion goes fast to a region 
where the field $v / \sqrt{\log(T)}$ has a
fixed deterministic profile, and stays there during this time interval.
 To each scenario corresponds
a lower bound of the type
\[
\lim_{\epsilon \rightarrow 0}
\liminf_{T \rightarrow \infty}
\frac{1}{T} \log P_0 \cro{|Y_T-y| \leq \epsilon} 
\geq - \II_n(y) \, .
\]
The family of functions $\II_n$ is decreasing, and 
the limit $\II(y) \triangleq \lim_{n \rightarrow \infty} \II_n(y)$
is convex. This enables us to identify $\II$ and the upper bound $\JJ$.

The paper is organized as follows. In section \ref{cadre}, we
introduce the notations and state the main result. In
section \ref{bornesup}, we prove the large deviations upper bound.
In section \ref{borneinf}, we establish the  large deviations
lower bound. Finally, section \ref{rate} investigates the
link between the decay of correlation, and the behavior
 of the rate function near the origin.

As a concluding remark, we would like to say that the paper
is written for a diffusion in $\R^2$, but that with a little more work,
all could be written in higher dimensions, 
as soon as the shear flow structure is preserved.

\vspace{1cm}
\section{Notations and results.}
\label{cadre}
\indent
   In all the sequel, when $I$ is a domain of $\R$, $\MM(I)$,
and $\MM_1(I)$ will denote respectively the set of finite measures on $I$,
and the sets of probability measures on $I$. $\CC^{\infty}_c(I)$,
$\CC(I)$, $H^1_0(I)$ will be respectively the set of infinitely
differentiable functions with compact support in $I$, the set
of continuous functions, and the Sobolev space obtained by completion
of $\CC^{\infty}_c(I)$ under the norm $\nor{f}_{H^1_0(I)}
= \int_I f^2(x) \, dx + \int_I (f')^2(x) \, dx$. Finally, for
all $p \in [1, \infty]$, $\nor{f}_p$ will denote the norm 
of the function $f$ in $L^p(I)$.

  Let $(v(x), x \in \R)$ be a centered stationary Gaussian
field with values in $\R$, defined on a probability
space $(\XX, \GG, \nu)$. Brackets  will denote the expectation with respect
to $\nu$, so that the covariance function of $v$ is defined  by
$K(x-y) \triangleq  <v(x) v(y)>$. 

Let $(B_t; t \in [0,1])$ be a standard 
Brownian motion defined on a probability
space $(\Omega, \AA, P)$. Expectation with respect to $P$ is denoted by 
$E$.

Our main result is  a full  large deviations principle 
for the random variable
$$Y_T \triangleq  \frac{1}{T \sqrt{\log(T)}} 
\int_{0}^{T} v\pare{B_s}  \, ds \,\,
.$$

Before stating the result, we introduce some 
assumptions and recall some standard results about the
Gaussian field $v$.

\subsection{The Gaussian field.}
\label{champ}
We assume that $v$ has a spectral density $h$ 
such that for some $\alpha >0$,
\begin{equation} 
\label{densiteL1}
\int_{\R} 
(1+|\lambda|^{\alpha}) \, 
h(\lambda) \, d\lambda < +\infty \, .
\end{equation}
Then, the covariance $K(x) = \int_{\R} e^{i \lambda x} \, h(\lambda)
\, d\lambda$
is a continuous function on $\R$, which attains its
maximal value at $0$. 
Moreover, $K(x) \rightarrow 0$ when $|x| \rightarrow  \infty$,
and $K$ is H\"older
 continuous of order $\alpha$, so that $v$ has a version which is
 $\beta$-H\"older continuous for $0<\beta<\frac{\alpha}{2} $.
Moreover, as it is well known for Gaussian fields,
\begin{equation}
\label{max-gaus}
\nu-\mbox{a.s.,} \,\,\, \limsup_{L \rightarrow +\infty} 
\frac{\Max \acc{|v(x)|: x\in [-L,L]} }{\sqrt{2K(0) \log(L)}} \leq 1 \, .
\end{equation}

We present now a splitting of $v$ into the sum of two Gaussian 
stationary processes, one of which having finite correlation
length. 
This splitting is constructed in \cite{gartner-konig-molchanov},
and goes as follows.

Let $g$ be the $L^2$-Fourier transform of $\sqrt{h}$. We can assume
that $v(x) = \int_{\R} g(x-y) \, dZ(y)$, where $Z$ is a Brownian
motion on $\R$.  Let $\psi: \R \mapsto[0,1]$ 
be a smooth even function, such that $\psi=0$ outside
$]-1/2;1/2[$, and $\psi=1$ on $[-1/4;1/4]$. Let $\psi_L(x) \triangleq 
\psi(\frac{x}{L})$, $g_L(x) \triangleq \psi_L(x) g(x)$,
and $\tilde{g}_L(x) \triangleq g(x)-g_L(x)$. 
This splitting of $g$ yields a corresponding
splitting of $v=v_L + \tilde{v}_L$, where
\begin{equation}   
\label{dec-champ}
 v_L(x) = \int g_L(x-y) \, dZ(y) \,\, , \,\,\, 
\tilde{v}_L(x) = \int \tilde{g}_L(x-y) \, dZ(y) \,\, .  
\end{equation}  
$v_L$ and $\tilde{v}_L$ are clearly stationary Gaussian processes.
The support of $K_L(x) = \bra{v_L(x) v_L(0)} = g_L \star g_L(x)$ (where
$\star$ denotes the convolution operator), 
is included in $[-L;L]$. 

Note also that
$\tilde{K}_L(0)=\bra{\tilde{v}_L(0) \tilde{v}_L(0)} 
= \int (1-\psi_L(x))^2 g^2(x) \, dx$ tends to $0$ when $L$ goes to infinity.

Moreover, if $\check{f}$ denotes the inverse
Fourier transform of $f$, 
\[ 
\begin{array}{l}
\int |\lambda|^{\alpha} |\check{g_L}(\lambda)|^2 \, d \lambda  
\\ \hspace{1cm}
 = \int |\lambda|^{\alpha} |\check{\psi_L} \star \sqrt{h}|^2
	 \, d \lambda 
\\ \hspace{1cm}
= \int d \lambda_1 \, d \lambda_2 \, \check{\psi_L}(\lambda_1) 
			\check{\psi_L}(\lambda_2)
   \int d \lambda \, |\lambda|^{\alpha} \sqrt{h}(\lambda-\lambda_1)
	 \sqrt{h}(\lambda - \lambda_2) 
\end{array}
\]
But 
\[
\begin{array}{ll}
\int d \lambda \, |\lambda|^{\alpha} \sqrt{h}(\lambda-\lambda_1)
	 \sqrt{h}(\lambda - \lambda_2) 
& \leq \prod_{i=1,2}
\pare{\int d \lambda \, |\lambda|^{\alpha} h(\lambda-\lambda_i)}^{1/2}
\\
& \leq C  \prod_{i=1,2}	\pare{ |\lambda_i|^{\alpha} +
	\int d \lambda \, |\lambda|^{\alpha} h(\lambda) }^{1/2} \, ,
\end{array} 
\]
so that 
\[
\int |\lambda|^{\alpha} |\check{g_L}(\lambda)|^2 \, d \lambda 
 \leq C \pare{ \int d \lambda \, (1 + |\lambda|^{\alpha})^{1/2} 
	|\check{\psi_L}(\lambda)|}^2 
<  \infty \, ,
\]
since   $\check{\psi_L}$ decreases faster than any polynomial at infinity.
Thus, $v_L$ has a H\"older continuous version, and so does $\tilde{v}_L$.

\subsection{The large deviations principle.}
Let us now define the rate function $\JJ$ appearing in the large deviations 
principle.
  When $f$  is a function of the Sobolev space $H^1(\R)$,
$K \star f^2$ is the continuous
function obtained by convolution of the covariance kernel $K$
and $f^2$, so that 
$$\pare{K \star f^2, f^2} = \int_{\R^2} K(x-y) f^2(x) f^2(y) \, dx \, dy
= \int_{\R} |\widehat{f^2}(\lambda)|^2 h(\lambda) \, d\lambda  \,\, .
$$
For any $\alpha \in \R$, let
\begin{equation}
\label{loglaplace}
\Lambda(\alpha) \triangleq 
\Sup 
\acc{ |\alpha| \sqrt{2 (K \star f^2,f^2)} - \frac{1}{2} \nor{f'}_2^2:
f \in H^1(\R), \nor{f}_2=1} \, , 
\end{equation}
and for any $y \in \R$
\begin{equation}
\label{suprate}
\JJ(y) \triangleq 
\Sup
\acc{\alpha y - \Lambda(\alpha):  \alpha \in \R} \, .
\end{equation}

We are now able to state the main result of the paper.

\begin{theo}	
\label{ldp-theo} 
Assume \refeq{densiteL1}. Then, $\nu$-a.s, for any measurable subset 
$E$ of $\R$,
\begin{equation}
\label{uby-eq}
\limsup_{T \rightarrow \infty} \frac{1}{T} 
\log P_0 \cro{Y_T  \in E}
\leq - \inf_{y \in \bar{E}} \, \JJ(y) \,\, , 
\end{equation}
\begin{equation}
\label{lby-eq}
\liminf_{T \rightarrow \infty} \frac{1}{T} 
\log P_0 \cro{Y_T  \in E}
\geq - \begin{array}[t]{c} \inf \\[-10pt]
	 \scriptstyle{y} \in \hspace{-2mm} \inte{\scriptstyle{E}} 
	\end{array}  \, \JJ(y) \,\, .
\end{equation}
$\JJ$ is even, convex, and lower semicontinuous. 
$\JJ(y) < \infty$ for $|y| < \sqrt{2 K(0)}$, and $\JJ(y) = +\infty$
for $|y| > \sqrt{2 K(0)}$. Moreover, $\JJ(0)=0$, 
and  $\JJ$ is increasing on $\R^+$.
\end{theo}

As a corollary of the large deviations for $Y$, we obtain the large
deviations for $X_2$ with the same rate function.

\begin{cor} 
\label{ldpx-cor} 
Assume \refeq{densiteL1}. Then,
the estimates \refeq{uby-eq} and \refeq{lby-eq} hold
when $X_{2,T}/(T \sqrt{\log(T)})$ replaces $Y_T$.
\end{cor}

We provide some more informations on $\JJ$, relating the
decay of correlation of the field $v$, and the behavior near
the origin of $\JJ$.
 
\begin{prop}	
\label{prop-J}.
\begin{enumerate} 
\item  Assume that for some $\beta \in ]0,1[$,
$\sous{\limsup}{|x| \rightarrow \infty} |x|^{\beta}|K(x)|
< \infty$, then 
\[
\liminf_{y \rightarrow 0} \frac{\JJ(y)}{|y|^{4/\beta}} > 0 \, .
\]
Assume that $K \geq 0$ and 
$\sous{\liminf}{|x| \rightarrow \infty} |x|^{\beta} K(x) > 0$ 
(for some $\beta \in ]0,1[$), or that
$\limi{|x|\rightarrow \infty} |x|^{\beta} K(x) > 0$,
 then 
\[
\limsup_{y \rightarrow 0} \frac{\JJ(y)}{|y|^{4/\beta}}
< \infty \, .
\]
\item Assume that for some $\beta > 1$, 
$\limisup{|x| \rightarrow \infty} |x|^{\beta} |K(x)|  < \infty$, and
$\int K(x) \, dx \neq 0$. Then, $\limi{y \rightarrow 0} \frac{\JJ(y)}{y^4}$
exists in $]0,+\infty[$ .
\end{enumerate}
\end{prop}

\vspace{1cm}
\section{Proof of the upper bound.}
\label{bornesup}

The aim of this section is to prove \refeq{uby-eq}, and the same estimate
for $X_{2,T}$. We begin with the proof of the
properties of $\JJ$ stated in theorem \ref{ldp-theo}. 

\subsection{Proof of the properties of $\JJ$.}
$\JJ$ is convex and
l.s.c as the supremum of affine  functions. $\JJ$ is even because 
$\Lambda$ is even. We restrict therefore the study of $\JJ$ to 
$\R^+$. For $y \in \R^+$,
$$ \sous{\Sup}{\alpha \leq 0} (\alpha y - \Lambda(\alpha))
 = \sous{\Sup}{\alpha \geq 0} (- \alpha y - \Lambda(\alpha))
 \leq \sous{\Sup}{\alpha \geq 0} (\alpha y - \Lambda(\alpha)) \,\, ,
 $$
so that $\forall y \in \R^+$, 
$\JJ(y) = \Sup \acc{\alpha y - \Lambda(\alpha),\alpha \geq 0}$. 
The monotony of $\JJ$ is thus obvious.

Let us prove now that $\JJ(y) < \infty$ for $|y| < \sqrt{2 K(0)}$,
and $\JJ(y)=+\infty$ for $|y| > \sqrt{2 K(0)}$. For this purpose, note
that 
\begin{equation}
\label{K0=sup}
\Sup \acc{(K \star f^2, f^2): f \in H^1(\R), \nor{f}_2=1}  = K(0) \, .
\end{equation}
Indeed, on one hand,  $\forall f \in H^1(\R)$ 
such that $\nor{f}_2=1$,  
$(K \star f^2, f^2) \leq K(0)$.  On the other hand, let $f_0$ be
any function in $H^1(\R)$, such that $\nor{f_0}_2=1$, and let 
$\lambda >0$. $f_{0,\lambda}(x)=\sqrt{\lambda} f_0(\lambda x)$ 
is then a function in $H^1(\R)$, such
that $\nor{f_{0,\lambda}}_2=1$. Therefore,
\[
\begin{array}{ll}
\Sup \acc{ (K \star f^2,f^2): f \in H^1(\R), \nor{f}_2=1}
& 
\geq (K*f_{0,\lambda}^2, f_{0,\lambda}^2)  \\
& 
= \int K \pare{\frac{x-y}{\lambda}} f_0^2(x) f_0^2(y) \, dx \, dy \, ,
\end{array} 
\]
and  \refeq{K0=sup}  follows by letting $\lambda \rightarrow \infty$,
and dominated convergence.
Thus, $\Lambda(\alpha) \leq |\alpha| \sqrt{2 K(0)}$, and 
$\forall y > \sqrt{2 K(0)}$, 

\[
\JJ(y) \geq \sous{\Sup}{\alpha \geq 0}
\acc{\alpha (y-\sqrt{2 K(0)})  } = + \infty \,\,.
\] 

On the other side, 
for $0 \leq y < \sqrt{2 K(0)}$, \refeq{K0=sup} allows one 
to find  $f_y$ in $H^1(\R)$  such
that $\nor{f_y}_2=1$, and $y < \sqrt{ 2 (K \star f^2_y,f^2_y)}$.
We get then that 
\[
\JJ(y) \leq \sous{\Sup}{\alpha \geq 0} \acc{ \alpha y
- \alpha \sqrt{ 2 (K \star f^2_y, f_y^2)} + \Frac{1}{2} \nor{f'_y}_2^2}
= \Frac{1}{2} \nor{f'_y}_2^2  <+\infty \,\, .
\]

Let us now compute 
$\JJ(0)$. Since $\Lambda$ is
even and increasing on $\R^+$, 
\[
\JJ(0)  = - \Inf \{ \Lambda(\alpha); \alpha \in \R\}
= - \Lambda(0) = 0 \, .
\]

\subsection{Large deviations upper bound for $X_{2,T}$.}
We are going to prove that \refeq{uby-eq} implies the 
same estimate for $X_2$.
Let us then assume that \refeq{uby-eq} holds.
 Let $\delta  > 0$, and let 
$F^{\delta }=\{ y: \exists x \in F, |x-y| \leq \delta \}$.
$$P_0 \cro{\frac{X_{2,T}}{T \sqrt{\log(T)}} \in F}
\leq  P_0 \cro{Y_T \in F^{\delta }}
+ P_0 \cro{ \va{\frac{W_{2,T}}{T \sqrt{\log(T)}}} \geq \delta } \,\, .
$$
But $\lim_{T \rightarrow \infty} \frac{1}{T} 
\log P_0 \cro{ \va{\frac{W_{2,T}}{T \sqrt{\log(T)}}} \geq \delta } 
= - \infty$.
Therefore, \refeq{uby-eq} yields  that $\nu$-a.s., for
all closed subset $F$, and all $\delta > 0$,
$$\limsup_{T \rightarrow \infty} \frac{1}{T} \log
 P_0 \cro{\frac{X_{2,T}}{T \sqrt{\log(T)}} \in F} \leq 
- \sous{\Inf}{y \in F^{\delta}} \JJ(y) \,\,.
$$
The result follows from the goodness
 of the rate function $\JJ$, letting $\delta$ go to $0$.  

\subsection{Large deviations upper bound for $Y_T$.}
We prove now \refeq{uby-eq} in theorem \ref{ldp-theo}.

\vspace{.5cm}
\noindent 
{\bf Step 1. Restriction of the problem in a domain of size $T$.} \\
For $R>0$, let $I_{RT}$ be the interval $]-RT;+RT[$, and let $\tau_{RT}$
be the first time Brownian $B$ exits  $I_{RT}$.
 
\begin{Lemme} 
\label{step1-lem}
$\nu$-.a.s, for all measurable set  $F$ and all $R>0$, 
\begin{equation}
\begin{array}{l}

\limsup_{T \rightarrow \infty} \frac{1}{T}
\log P_0 \cro{ Y_T \in F}
\\ \hspace{2cm}
\leq \Max \acc{ \limsup_{T \rightarrow \infty} \frac{1}{T}
\log P_0 \cro{ Y_T \in F; \, \tau_{RT} > T}, - \frac{R^2}{2}} \,\, .
\end{array}
\end{equation}
\end{Lemme}

\noindent 
{\bf Proof.} 
\begin{equation}
P_0( Y_T  \in F) \leq P_0 \cro{Y_T \in F; \, \tau_{RT} > T}
+ P_0 \cro{ \Sup_{[0,T]} |B_s| \geq RT} \,\, .
\end{equation}
The well known estimate 
 $\limsup_{T \rightarrow \infty} \frac{1}{T} 
\log P_0 \cro{ \Sup_{[0,T]} |B_s| \geq RT} \leq - \frac{R^2}{2}$
yields the  result.
\hspace*{\fill} \rule{2mm}{2mm} 

\vspace{1cm}
\noindent
{\bf Step 2. Spectral estimates of Schr\"odinger semigroups.} \\
To prove the upper bound, we use the
G\"artner-Ellis method, and  we have to study the large time
asymptotic of 
\begin{equation}
\Lambda_{T}\pare{\frac{\alpha v}{\sqrt{\log(T)}},I_{RT}}
= E_0 \cro{
\exp \pare{ \int_0^T  \frac{\alpha v(B_s)}{\sqrt{\log(T)}} \, ds};\,
\tau_{RT} > T} \,\, .
\end{equation}
It is well known that this reduces to study the principal eigenvalue  
of the random operator 
$\LL(f) = \frac{1}{2} f''+ \alpha \frac{v}{\sqrt{\log(T)}} f$,
with Dirichlet conditions on the boundary of $I_{RT}$. 

 In all the sequel, when  $D$ is a  bounded domain of $\R$, 
and $V: D \mapsto \R$ is a 
bounded measurable function,  
 we will write $\lambda(V,D)$ for the principal
eigenvalue of the operator $1/2 \triangle + V$, with Dirichlet
boundary condition on  $D$.
\[
\lambda(V,D) 
 \triangleq 
\inf 
 \acc{ \frac{1}{2} \int_D (f')^2(x) \, dx 
- \int_D V(x) f^2(x) \, dx: 
f \in \CC^{\infty}_c(D), \nor{f}_2=1}  
\]
Since any sequence $(f_n)$ which is bounded in $H^1_0(D)$
 has a subsequence which converges strongly in $L^2(D)$
and weakly in  $H^1_0(D)$, one also has
\[
\lambda(V,D) = 
\Min 
 \acc{ \frac{1}{2} \int_D (f')^2(x) \, dx 
- \int_D V(x) f^2(x) \, dx :
f \in H^1_0(D), \nor{f}_2=1}
\]

\vspace{.5cm}
\noindent 
 In these notations, the task at hand is
 to study the behavior for large $T$ of
$\lambda (\alpha v/\sqrt{\log(T)} , I_{RT})$. 
To this end, we recall proposition 1 of \cite{gartner-konig}, which
compares this eigenvalue, with the minimum of the principal eigenvalues
in balls of fixed size.
  
\begin{Lemme}
\label{GK-lem}
(Proposition 1 of \cite{gartner-konig}).\\
$\forall r \geq 2$, there exists a continuous $2r$-periodic
function $\Phi_r : \R \mapsto \R^+$, with support included
in $\cup_{k \in \Z} ((2k+1)r + ]-1;1[)$, such that
for all $R>r$, for all H\"older continuous $V: \R \mapsto \R$,
for all $\theta \in I_{2r}$
\begin{equation}
\label{GK-eq}
\lambda(V-\Phi^{\theta}_r,I_R) \geq
\Min \acc{ \lambda(V, z+I_{2r+1}): z \in (2r \Z) \cap I_{R+r}}\,\,,
\end{equation}
where $\Phi^{\theta}_r(x) = \Phi_r(x-\theta)$.\\
Moreover, $\frac{1}{|I_r|} \int_{I_r} \Phi_r(x) \, dx \leq \frac{K}{r}$,
where the constant $K$ is independent of $r$.
\end{Lemme}

\vspace{.5cm}
\noindent 
We deduce from this the following lemma.

\begin{Lemme}
\label{Lambda-leq-minvp}
There exists a constant $K$ such that $\nu$-a.s., for all $r \geq 2$,
 $\forall \alpha \in \R$, $\forall R >0$ 
\[
\begin{array}{l}
\limisup{T \rightarrow \infty} \, 
\frac{1}{T} \log  \Lambda_{T} \pare{\frac{\alpha v}{\sqrt{\log(T)}},I_{RT}}
\\ \hspace*{1cm}
\leq \Frac{K}{r}  
- \limiinf{T \rightarrow \infty} 
\Min \acc{ 
\lambda \pare{\frac{\alpha v}{\sqrt{\log(T)}},z+I_{2r+1}}:
z \in (2r \Z) \cap I_{RT+r}}  \,\, .
\end{array} 
\]
\end{Lemme}

\noindent 
{\bf Proof.} 
We use  the same trick as in \cite{gartner-konig} and 
\cite{gartner-konig-molchanov}. Let $\Phi_r$ be
the function introduced in lemma \ref{GK-lem}. By periodicity of
$\Phi_r$, $\frac{1}{|I_r|} \int_{I_r} \Phi_r(\theta + B_s) \, d\theta 
= \frac{1}{|I_r|} \int_{I_r} \Phi_r(\theta) \, d\theta  \leq \frac{K}{r}$.
By Jensen inequality, we obtain then that 
\[
\Lambda_{T}\pare{\frac{\alpha v}{\sqrt{\log(T)}},I_{RT}}
\leq \exp\pare{\frac{KT}{r}}
\, \, \frac{1}{|I_r|} 
\int_{I_r} d \theta \,\, 
\Lambda_T \pare{\frac{\alpha v}{\sqrt{\log(T)}} - \Phi^{\theta}_r,I_{RT}}
\,\, .
\]
We use then the usual bounds on Schr\"odinger  semigroups
 in terms of their principal eigenvalue (see
for instance theorem 1.2 in chapter 3 of \cite{sznitman}).  
\[
\begin{array}{l}
\Lambda_{T} \pare{\frac{\alpha v}{\sqrt{\log(T)}},I_{RT}} 
\\
 \leq C e^{\frac{KT}{r}}
\pare{
1 + \sqrt{ T \sous{\sup}{\theta \in I_r} 
\lambda (\scriptstyle{\frac{\alpha v}{\sqrt{\log(T)}}-\Phi^{\theta}_r, 
I_{RT}})}}
\exp \pare{- T \sous{\inf}{\theta \in I_{r}}
\lambda (\scriptstyle{\frac{\alpha v}{\sqrt{\log(T)}} -\Phi^{\theta}_r,
I_{RT}})} 
\\
\leq 
C e^{\frac{KT}{r}}
\pare{  
1 + \sqrt{T} 
( \nor{\Phi_r}_{\infty}^{\frac{1}{2}}  + 
	\frac{\max_{I_{RT}}|\alpha v|^{\frac{1}{2}}}{\log(T)^{\frac{1}{4}}})} 
\\
\hspace*{\fill} \exp \pare{- T 
\Min \acc{\lambda \pare{\frac{\alpha v}{\sqrt{\log(T)}},
z+I_{2r+1}}:
z \in 2r\Z \cap I_{RT+r}} }
\end{array} 
\]
The conclusion follows from \refeq{max-gaus} and \refeq{GK-eq}.
\hspace*{\fill} \rule{2mm}{2mm}

\vspace{1cm}
\noindent
{\bf Step 3. $\nu$-a.s. behavior of 
$\Min \acc{
\lambda \pare{\frac{\alpha v}{\sqrt{\log(T)}},z+I_{2r+1}}:
z \in 2r\Z \cap I_{RT+r}}$. }\\ 
This is done via a Borel-Cantelli argument. Using the
stationarity of $v$, note that the random variables
$\acc{\lambda \pare{\alpha v/\sqrt{\log(T)},z+I_{2r+1}}, 
z \in 2r\Z \cap I_{RT+r}}$ have the same law. The next lemma 
gives some estimates for this law.

\begin{Lemme} 
\label{queue-vp-lem}
Let $c \triangleq \Min \{\frac{1}{2} \int (f')^2 \, dx :
f \in H^1_0(I_1), \int f^2 =1\}$. Let us define for 
all $x \in \R$,  and  $r>0$
\begin{equation}
\label{jr-def}
J_r(x) \triangleq \left\{ \begin{array}{ll}
\inf \acc{ 
   \Frac{\pare{\frac{1}{2} \int(f')^2 -x}^2}
			{2 (K \star f^2,f^2)}:
 f \in H^1_0(I_r), \int f^2 =1} 
& \mbox { if } x < \frac{c}{r^2} \, , \\
0 & \mbox{ otherwise.}
\end{array}
\right.
\end{equation}
Then, $\forall r>0$,  $\forall x \in \R$, 
\begin{equation}
\label{queue-vp-eq}
\lim_{T \rightarrow \infty} \Frac{1}{\log(T)}
\log \nu \cro{\lambda \pare{\frac{\alpha v}{\sqrt{\log(T)}},I_r} 
\leq x} = -\frac{J_r(x)}{\alpha^2} \,\, .
\end{equation}
\end{Lemme}

\noindent
{\bf Proof.} Let $f$ be any function in $H^1_0(I_r)$ such
that $\int f^2=1$. Then
\[
\nu \cro{\lambda \pare{\frac{\alpha v}{\sqrt{\log(T)}},I_r}  \leq x}  
 \geq \nu \cro{ (\alpha v,f^2) \geq 
	\sqrt{\log(T)} (\frac{1}{2} \nor{f'}_2^2  - x)}  \,\, .
\]
But $(\alpha v,f^2) \sim \NN(0, \alpha^2 (K \star f^2, f^2))$, so
that
\[
\begin{array}{l}
\limiinf{T \rightarrow \infty}
\frac{1}{\log(T)} \log 
\nu \cro{\lambda\pare{\frac{\alpha v}{\sqrt{\log(T)}},I_r}  \leq x}   \\
\hspace*{2cm}  \geq
\left\{ \begin{array}{ll}
0 & \mbox{ for } x \geq \frac{1}{2} \nor{f'}^2_2 \, , \\
-   \Frac{ (\frac{1}{2} \nor{f'}_2^2 -x)^2}
		   {2 \alpha^2 (K \star f^2,f^2)} 
& \mbox{ for }  x <  \frac{1}{2} \nor{f'}^2_2 \, .
\end{array}
\right.
\end{array}
\]
Taking the supremum over all functions $f \in H^1_0(I_r)$
such that $\nor{f}_2=1$, yields
\begin{equation}
\label{lb-fr-vp}
\liminf_{T\rightarrow \infty} \frac{1}{\log(T)}  
\log \nu \cro{\lambda\pare{\frac{\alpha v}{\sqrt{\log(T)}},I_r}  \leq x}  
 \geq - \frac{J_r(x)}{\alpha^2} \,\, .
\end{equation} 

\vspace{.5cm}
\noindent 
We are now going to prove the upper bound. 
To this end, note that $\lambda(\cdot,I_r): \CC(\bar{I_r}) 
\rightarrow \R$ 
is continuous (the topology in $\CC(\bar{I_r})$
being given by the supremum norm).
Indeed, first $\lambda(\cdot,I_r)$ is u.s.c as infimum 
of continuous functions. Secondly, we prove the lower semicontinuity:
 let then
$(v_n, n\in \N)$ be a sequence in $\CC(\bar{I_r})$ converging to $v$.
For all $n \in \N$, let $f_n$  realize  the infimum  in $\lambda(v_n,I_r)$.
Since $\lambda(v_n,I_r) \leq - \Min_{Ir} v_n$, and $\nor{v_n-v}_{\infty} 
\rightarrow 0$, the sequence $(f_n)$ is bounded in $H^1_0(I_r)$, and
admits therefore a subsequence converging strongly in $L^2(I_r)$
and weakly in  $H^1_0(I_r)$ to a function $f \in H^1_0(I_r)$. One 
obtains then that $\nor{f}_2 = \lim \nor{f_n}_2 =1$, 
$\liminf \nor{f'_n}_2 \geq \nor{f'}_2$, and $\lim (v_n,f_n^2)
= (v,f^2)$, so that  $\liminf_{n \rightarrow \infty} \lambda(v_n,I_r)
\geq \frac{1}{2} \nor{f'}^2_2 - (v,f^2) \geq \lambda(v,I_r)$.

Therefore, for all $r >0$, $x, \alpha \in \R$, 
$
F_r^{\alpha,x} \triangleq 
\acc{ u \in \CC(\bar{I_r}),  \lambda(\alpha u,I_r) \leq x}
$
is a closed subset of $\CC(\bar{I_r})$, and 
\[ 
\nu \cro{ \lambda\pare{\frac{\alpha v}{\sqrt{\log(T)}},I_r} \leq x}
= \nu \cro{ \frac{v}{\sqrt{\log(T)}} \in F^{\alpha,x}_r} \,\, .
\]
We now use the large deviations in $\CC(\bar{I_r})$ of the Gaussian field
$ v/\sqrt{\log(T)}$ to deduce that 
\[
\begin{array}{l}
\limsup_{T \rightarrow \infty} \Frac{1}{\log(T)} 
\, \nu  \cro{ \lambda\pare{\frac{\alpha v}{\sqrt{\log(T)}},I_r} \leq x}
\\
\hspace{2cm}
\leq -  \inf \acc{ K_r^*(u): u \in \CC(\bar{I}_r),
\lambda(\alpha u,I_r) \leq x } \,,
\end{array}
\]
where 
\begin{equation}
\label{kstar-def}
K^*_r(u) \triangleq 
\sup
\acc{ (u,\mu) - \frac{1}{2} (K \star \mu,\mu):
\mu \in \MM(I_r)} \,\, .
\end{equation}
Note that 
\begin{equation}
\label{k*carre}
\begin{array}{ll}
K^*_r(u) & = \sous{\Sup}{\mu \in \MM(I_r)} 
		\sous{\Sup}{ m \in \R }
  		\acc{ m(u,\mu) - \Frac{m^2}{2} (K  \star \mu,\mu)} \\
& =  \sous{\Sup}{\mu \in \MM(I_r)} 
               \acc{ \Frac{(u,\mu)^2}{2 (K \star \mu,\mu)}}
(\mbox{with the convention } \frac{0}{0}=0).
\end{array}
\end{equation}	
Hence $\forall \alpha \in \R$, $K^*_r(\alpha u)=\alpha^2 K^*_r(u)$,
and 
\[
\begin{array}{l}
\limsup_{T \rightarrow \infty} \Frac{1}{\log(T)} 
\, \nu  \cro{ \lambda\pare{\frac{\alpha v}{\sqrt{\log(T)}},I_r} \leq x}
\\
\hspace{2cm}
\leq - \frac{1}{\alpha^2} \inf \acc{ K_r^*(u): u \in \CC(\bar{I}_r),
\lambda(u,I_r) \leq x } \,,
\end{array}
\]

It remains now to show that
\begin{equation} 
\label{bs=bi}
\inf \acc{ K_r^*(u):  u \in \CC(\bar{I_r}),
\lambda(u,I_r)  \leq x} \geq  J_r(x) \, \, .
\end{equation} 
We can restrict ourselves to the case where
$x < \frac{c}{r^2}$.  Let $u \in \CC(\bar{I_r})$ be such that 
$\lambda(u,I_r) \leq x$. Let $f_u \in H^1_0(I_r)$ be such that
$\nor{f_u}_2 =1$ and $\lambda(u,I_r) =  \frac{1}{2} \nor{f'_u}^2_2
- (u,f_u^2)$.
It follows from \refeq{k*carre} that 
$K_r^*(u) \geq \Frac{(u,f_u^2)^2}{2 (K \star f^2_u,f^2_u)}$.
But 
\[
(u,f^2_u) = - \lambda(u,I_r) + \frac{1}{2} \nor{f'_u}^2_2 
\geq \frac{1}{2} \nor{f'_u}^2_2 - x \,\, .
\]
Moreover, $x < \frac{c}{r^2} \leq \frac{1}{2} \nor{f'_u}^2_2$ by
definition of the constant $c$. Thus,
\[ 
K_r^*(u) \geq \Frac{\pare{\frac{1}{2} \nor{f'_u}_2^2 -x}^2}
{2 (K \star f^2_u, f^2_u)} \geq J_r(x) \,\, .
\]
Taking the infimum over functions $u$ such that $\lambda(u,I_r) \leq x$
yields then \refeq{bs=bi}.
\hspace*{\fill} \rule{2mm}{2mm}

\vspace{.5cm}
Lemma \ref{queue-vp-lem} allows one to prove

\begin{Lemme}
\label{comp-min-vp-lem}
$\forall \alpha \in \R$, and  $\forall r \geq 2$,  let 
\begin{equation}
\label{lambdar-def}
\Lambda(\alpha, r) \triangleq  
\Sup \acc{ |\alpha| \sqrt{2(K\star f^2,f^2)} - \frac{1}{2} \nor{f'}_2^2:
  f \in H^1_0(I_{2r+1}) ,\, \nor{f}_2 =1} \, .
\end{equation}
Then, $\forall \alpha \in \R$, $\forall R>0$ and  $\forall r \geq 2$,
$\nu$-a.s.,
\begin{equation}
\label{comp-min-vp-eq} 
\liminf_{T \rightarrow \infty} \,
\Min\acc{
\lambda \pare{\frac{\alpha v}{\sqrt{\log(T)}},z + I_{2r+1}}:
z \in (2r\Z)\cap I_{RT+r} }
\geq - \Lambda(\alpha,r) \, .
\end{equation}
\end{Lemme}
 
\noindent 
{\bf Proof.} We use Borel-Cantelli lemma. We assume that 
$\Lambda(\alpha, r) < \infty$, otherwise there is nothing to prove. Let
$\epsilon > 0$ be fixed.
\[
\begin{array}{l}
\nu \cro{\Min \acc{
\lambda \pare{\frac{\alpha v}{\sqrt{\log(T)}},z + I_{2r+1}}:
z \in (2r\Z)\cap I_{RT+r} }
\leq - \Lambda(\alpha,r) - \epsilon }
\\
\hspace*{1cm}
\leq \sum_{ z \in (2r\Z)\cap I_{RT+r} }
\hspace{-.3cm}
	\nu \cro{\lambda\pare{\frac{\alpha v}{\sqrt{\log(T)}}
		,z + I_{2r+1}}
\leq - \Lambda(\alpha,r) - \epsilon } 
\\
\hspace*{1cm}
\leq C(1 + \frac{RT}{r})  
\nu \cro{\lambda\pare{\frac{\alpha v}{\sqrt{\log(T)}},I_{2r+1}}
\leq - \Lambda(\alpha,r) - \epsilon }  \mbox{ by stationarity. }
\end{array}
\]
Thus, by lemma \ref{queue-vp-lem} ,
\begin{equation}
\label{queue-min}
\begin{array}{l}
\limisup{T \rightarrow \infty} \frac{1}{\log(T)}
\log \nu \cro{\sous{\Min}{z \in (2r\Z)\cap I_{RT+r}}
\hspace{-.3cm}
 \lambda \pare{\frac{\alpha v}{\sqrt{\log(T)}},z + I_{2r+1}}
\leq - \Lambda(\alpha,r) - \epsilon }
\\
\hspace*{2cm}
\leq 1- \Frac{J_{2r+1}(-\Lambda(\alpha,r) - \epsilon)}{\alpha^2}  \,\, .
\end{array}
\end{equation}

We claim that 
\begin{equation}	
\label{inverse-de-J}
x < -\Lambda(\alpha,r) \Leftrightarrow  J_{2r+1}(x) > \alpha^2 \, .
\end{equation}
The only point to note in order to prove \refeq{inverse-de-J}
is that the infimum in \refeq{jr-def}, and the
supremum in \refeq{lambdar-def} are actually reached, since again
any majorizing sequence will be bounded in $H^1_0(I_{2r+1})$,
and $f \in L^2(I_{2r+1}) \mapsto (K \star f^2,f^2)$ is
continuous.
Hence,
\[
\begin{array}{l}
x < -\Lambda(\alpha,r)
\\ \hspace{1cm}
 \Leftrightarrow \forall f \in H^1_0(I_{2r+1}), \nor{f}_2=1, \,\,\,
|\alpha| \sqrt{ 2 (K \star f^2,f^2)} < \Frac{1}{2} \nor{f'}_2^2 - x 
\\ \hspace{1cm}
 \Leftrightarrow 
 	 \forall f \in H^1_0(I_{2r+1}), \nor{f}_2=1, 
	 \left\{ \begin{array}{l}
		 \Frac{1}{2} \nor{f'}_2^2 - x > 0 \\
		\Frac{(\frac{1}{2} \nor{f'}_2^2 - x)^2}{2 (K \star f^2,f^2)}
		> \alpha^2 
	\end{array} \right.
\\ \hspace{1cm} 
\Leftrightarrow J_r(x) > \alpha^2 \,.
\end{array}
\]

It follows then from \refeq{queue-min}, \refeq{inverse-de-J},
and Borel-Cantelli lemma applied along the sequence $T_n = 2^n  $,
that $\forall \alpha \in \R$, $\forall r \geq 2$, $\forall R > 0$,
$\nu$-a.s.,
\[
\liminf_{n \rightarrow \infty} \,
\Min \acc{
\lambda \pare{\frac{\alpha v}{\sqrt{\log(T_n)}},z + I_{2r+1}}:
z \in (2r\Z)\cap I_{RT_n+r} }
\geq - \Lambda(\alpha,r) \,\, ,
\]
To end the proof of lemma \ref{comp-min-vp-lem}, note that
for $T$ sufficiently large, and $n$ such that $T_n \leq T < T_{n+1}$,
\[
\begin{array}{l}
\Min \acc{
\lambda\pare{\frac{\alpha v}{\sqrt{\log(T)}},z + I_{2r+1}}:
z \in (2r\Z)\cap I_{RT+r} }
\\
\geq 
\Min \acc{
\lambda \pare{\frac{\alpha v}{\sqrt{\log(T_{n+1})}},z + I_{2r+1}}:
z \in (2r\Z)\cap I_{RT_{n+1}+r} }
- \frac{\sous{\max}{I_{RT_{n+1}+3r+1}} 
\hspace{-.3cm}
|\alpha v|}{\log(2) n (n+1)} \,\, .
\end{array}
\]
The last term is $\nu$-a.s. of order $1/n$ by \refeq{max-gaus}.
\hspace*{\fill} \rule{2mm}{2mm}

\vspace{.5cm}

Concerning lemma \ref{comp-min-vp-lem}, we would like to
underline that  using the decorrelation properties of the
field $v$, and  Borel Cantelli inverse lemma, it is possible
to prove that $-\Lambda(\alpha,r)$ is in fact the a.s. limit
when $T \rightarrow \infty$ of 
$\Min \acc{
\lambda\pare{(\alpha v)/\sqrt{\log(T)},z+I_{2r+1}}:
z \in (2r\Z)\cap I_{RT+r}}
$.

\vspace{.5cm}
At this point, putting lemma \ref{Lambda-leq-minvp} 
and lemma \ref{comp-min-vp-lem} together,
we have proved that there exists $K >0$ such that:
$\forall r \geq 2$, $\forall R >0$, $\forall \alpha \in \R$,
$\nu$-a.s.,
\[
\limsup_{T \rightarrow \infty } \Frac{1}{T}
\log \Lambda_T \pare{\frac{\alpha v}{\sqrt{\log(T)}},I_{RT}}
\leq \Frac{K}{r} + \Lambda(\alpha,r) \, .
\]
Taking the limit $r \rightarrow \infty$ along subsequences,
we obtain that $\nu$-a.s., $\forall \alpha \in \Q$,
$\forall R \in \Q^+$,
\begin{equation}
\label{limsup-Lambda}
\limsup_{T \rightarrow \infty } \Frac{1}{T}
\log \Lambda_T\pare{\frac{\alpha v}{\sqrt{\log(T)}},I_{RT}}
\leq  \Lambda(\alpha) \, .
\end{equation}

\vspace{.5cm}
\noindent
{\bf Step 4. Conclusion.}\\ 
 It is now routine to obtain from
\refeq{limsup-Lambda} the weak large deviations
upper bound (i.e. the upper bound for compact sets). 
\refeq{uby-eq} follows then from the exponential tightness of $Y$ (lemma
\ref{tension}).

\begin{Lemme} (weak large deviations upper bound). \\
$\nu$-a.s., $\forall y \in \R$, 
\[
\lim_{\epsilon \rightarrow 0} \, 
\limsup_{T \rightarrow \infty} \, 
\Frac{1}{T} \log
P_0 \cro{Y_T \in [y-\epsilon, y + \epsilon]}
\leq - \JJ(y) \, \, .
\]
\end{Lemme}

\noindent
{\bf Proof.} We treat only the case $y > 0$. By lemma \ref{step1-lem},
$\forall \epsilon < y$, $\forall \alpha >0$, and $\forall R>0$
\[
\begin{array}{l}
\limisup{T \rightarrow \infty} 
\Frac{1}{T} \log
P_0 \cro{Y_T \in [y-\epsilon, y + \epsilon]}
\\ \hspace{1cm} \leq
\Max \cro{\Frac{-R^2}{2},
\limisup{T \rightarrow \infty} 
\Frac{1}{T} \log P_0 \cro{Y_T \in [y-\epsilon, y + \epsilon]; 
		\tau_{RT} \geq T}}
\\ \hspace{1cm} \leq 
\Max \cro{\Frac{-R^2}{2},
-\alpha (y-\epsilon)  +
\limisup{T \rightarrow \infty} 
 \Frac{1}{T} \log E_0 \cro{ e^{\alpha T Y_T}; \tau_{RT} \geq T}}
\end{array}
\]
Therefore, $\nu$-a.s., $\forall y > 0$,
$\forall \epsilon < y$, $\forall R \in \Q^+$,
\[
\begin{array}{l}
\limsup_{T \rightarrow \infty}
\Frac{1}{T} \log
P_0 \cro{|Y_T -y | \leq \epsilon}
\\ \hspace{1cm}
\leq 
\Max \cro{\Frac{-R^2}{2}, 
- \sup \acc{\alpha (y-\epsilon) 
- \Lambda(\alpha):\alpha \in \Q^+} } \, .
\end{array}
\]
Note that by continuity of $\Lambda$, the supremum
on $\Q^+$, is a supremum on $\R^+$. Thus, \refeq{uby-eq} is
obtained by taking
the limit $R \rightarrow \infty$, then $\epsilon \rightarrow 0$, 
and by using the lower semi-continuity of $\JJ$.

\hspace*{\fill} \rule{2mm}{2mm}

\begin{Lemme} 
\label{tension} (exponential tightness).\\
$\nu$-a.s., $\forall L > \sqrt{2 K(0)}$,
\[
\limsup_{T \rightarrow \infty} \, 
\Frac{1}{T} \log
P_0 \cro{|Y_T| > L} 
\leq - \frac{L^2}{2} \, \, .
\]
\end{Lemme}

\noindent
{\bf Proof.}
Let $L > \sqrt{2 K(0)}$ be fixed. 
\[
P_0 \cro{|Y_T| > L}
\leq P_0 \cro{\tau_{LT} \leq T} + 
 \ind_{ \frac{\max_{I_{LT}}|v|}{\sqrt{\log(T)}} 
		> L } 
\leq C \exp(- \frac{L^2 T}{2}) + 
\ind_{ \frac{\max_{I_{LT}}|v|}{\sqrt{\log(T)}} > L } \,\, .
\]
By \refeq{max-gaus}, $\nu$-a.s., the indicator is null 
for $T$ sufficiently large. Therefore, $\forall  L > \sqrt{2 K(0)}$,
$\nu$-a.s., 
\[\limsup_{T \rightarrow \infty} \, 
\Frac{1}{T} \log
P_0 \cro{|Y_T| > L} 
\leq - \frac{L^2}{2} \, .
\]
Inverting the ``$\forall L$''  and the ``$\nu$-a.s'', is easily done
using the monotony of $L \mapsto P_0(|Y_T| > L)$. 
\hspace*{\fill} \rule{2mm}{2mm}


\vspace{1cm}
\section{Proof of the lower bound.}
\label{borneinf}

Here, we prove \refeq{lby-eq}, from
which the same assertion for $X_2$ is  easily deduced.

\subsection{a.s. behavior of the field with finite
correlation length.}

 As explained
in the introduction, the lower bound is obtained by forcing 
the Brownian motion to spend a certain amount of time in boxes
where the field $v/\sqrt{\log(T)}$ has a fixed profile. We need
therefore to describe the a.s.  behavior of this random
field. This is done in  the following lemma, assuming
that $K$ has compact support.  

\begin{Lemme}
\label{comp-ps-champ}
Assume that $K$ has compact support in $I_L$ for some $L > 0$. 
Let $\epsilon > 0$, $r > L$ and 
let $u$ be any function in $\CC(\bar{I_r})$ such that $K^*_r(u) < 1$.
Then $\nu$-a.s., for $T$ sufficiently large, 
$\exists z \in 2r\Z \cap I_{T/\log(T)}$ such that 
$\sous{Max}{y \in z+I_r} \va{ \frac{v(y)}{\sqrt{\log{T}}} - u(y-z)} 
\leq \epsilon$.
\end{Lemme}

\noindent
{\bf Proof.}
\[
\begin{array}{l}
\nu \cro{ \forall  z \in 2r\Z \cap I_{T/\log(T)} , \,
 \nor{ \frac{v(\cdot)}{\sqrt{\log{T}}} - u(\cdot - z)}_{\infty,z+Ir} \geq 
	\epsilon}
\\ \hspace{1cm}
\leq 
\nu \cro{ \forall  z \in 4r\Z \cap I_{T/\log(T)} , \,
 \nor{ \frac{v(\cdot)}{\sqrt{\log{T}}} - u(\cdot - z)}_{\infty,z+Ir} \geq 
	\epsilon} 
\end{array} 
\]
Since $K$ has compact support in $I_L$, and $r > L$, the random
variables $(\frac{v(z+I_r)}{\sqrt{\log{T}}},z \in 4r\Z \cap I_{T/\log(T)})$ 
are independent. Thus,
\[
\begin{array}{l}
\nu \cro{ \forall  z \in 2r\Z \cap I_{T/\log(T)} , \,
 \nor{ \frac{v(\cdot)}{\sqrt{\log{T}}} - u(\cdot - z)}_{\infty,z+Ir} \geq 
	\epsilon}
\\ \hspace{1cm}
\leq 
\pro{z \in 4r\Z \cap I_{T/\log(T)}}{}
\nu \cro{\nor{\frac{v(\cdot)}{\sqrt{\log{T}}}-u(\cdot-z)}_{\infty,z+Ir} \geq 
	\epsilon}
\\ \hspace{1cm}
\leq \pare{\nu \cro{\nor{\frac{v(\cdot)}{\sqrt{\log{T}}}-u}_{\infty,Ir} \geq 
	\epsilon}}^{2 \cro{\frac{T}{4r \log(T)}} +1}
\end{array}
\]
Let $\eta >0$ be such that $K^*_r(u) + \eta < 1$. Using 
the large deviations estimates of $v/\sqrt{\log(T)}$ , we obtain that for $T$ 
sufficiently large,
\[
\nu \cro{\nor{\frac{v(\cdot)}{\sqrt{\log{T}}}-u}_{\infty,Ir} \geq \epsilon} 
\leq 1 - T^{-K^*_r(u) - \eta} \,.
\]
Thus, for $T$ sufficiently large,
\[
\begin{array}{l}
\nu \cro{ \forall  z \in 2r\Z \cap I_{T/\log(T)} , \,
 \nor{ \frac{v(\cdot)}{\sqrt{\log{T}}} - u(\cdot - z)}_{\infty,z+Ir} \geq 
	\epsilon}
\\ \hspace{1cm}
\leq
\exp \pare{-  \pare{2 \cro{\frac{T}{4r \log(T)}} +1}T^{-K^*_r(u) - \eta}}
\sim \exp \pare{- \frac{T^{1-K^*_r(u)-\eta}}{2r \log(T)}} \, .
\end{array}
\]
The result follows by  Borel-Cantelli lemma applied 
along the sequence $T_n=n$.
\hspace*{\fill} \rule{2mm}{2mm} 

\subsection{Lower bounds for $Y_T$, with fixed profiles of the field.}

 From lemma \ref{comp-ps-champ}, we know that the field can be close
to $u$ with $K^*_r(u) < 1$, in a region a size $r$. Thus, for $n$ integer,
let 
\[
\UU(n,r) \triangleq \acc{ \vec{u} \in \CC(\bar{I}_r)^n; \Max_i K_r^*(u_i) <1}
\, ,
\]
be the $n$-tuples of admissible profiles.
  A lower bound for $P_0 \cro{ |Y_T - y| < \epsilon}$
is obtained by dividing $[0,T]$ into $n$ time intervals of length
$\alpha_i T$ ($0 \leq \alpha_i \leq 1, \sum_{i=1}^{n}  \alpha_i=1$). In
each time interval, we force the Brownian motion to go ``fast''
(say in a time of order $T/\log(T)$) from $I_1$ to a region
in $I_{T/\log(T)}$, where the field $v/\sqrt{\log(T)}$
is close to  $u_i$, to remain there
during $\alpha_i T - 2 T/\log{T}$, and then to return 
fast (in time of order $T/\log(T)$) to $I_1$.

Before stating the lower bound obtained in this way, we introduce
some notations. For any integer $n$, and any $r \in ]0,\infty]$,
define 
\[
\DD(n,r) \triangleq \acc{(\vec{\alpha},\vec{f} ) \in [0,1]^n \times 
H^1(I_r)^n:  \sum_{i=1}^{n} \alpha_i =1,  \, 
\forall \,  1 \leq i \leq n, \nor{f_i}_2 =1} \, ,
\]
\[
\DD(n) \triangleq  \DD(n,\infty), 
\mbox {and for } (\vec{\alpha},\vec{f}) \in \DD(n), \,
I_n(\vec{\alpha},\vec{f}) \triangleq
\frac{1}{2} \sum_{i=1}^{n} \alpha_i \nor{f_i'}_2^2
\,.
\]

\begin{Lemme}
\label{lb-u-lem}
Assume that $K$ has compact support in $I_L$. 
Then, $\forall r > L$, $\forall  \epsilon > 0$, $\forall n \in \N$,
$\forall \vec{u} \in \UU(n,r)$,
$\nu$-a.s., $\forall   y \in \R$, 
\begin{equation}
\label{lb-u}
\begin{array}{l}
\liminf_{T \rightarrow \infty} \frac{1}{T} \log P_0 
			\cro{ |Y_T - y| < \epsilon} 
\\ \hspace*{\fill}
\geq - 
\sous{\Inf}{z\in \BB(y,\epsilon)}
\sous{\Inf}{(\vec{\alpha},\vec{f}) \in \DD(n,r)}
\acc{I_n(\vec{\alpha},\vec{f}): 
\sum_{i=1}^{n} 	\alpha_i (u_i,f_i^2) = z }\, .
\end{array}
\end{equation}
\end{Lemme}

\noindent
{\bf Proof.} 
We begin with some more notations. 

For $0< S<T$, we will write $L_S^T$ for the occupation measure 
of $B$ between $S$ and $T$, 
$L_S^T = \frac{1}{T-S} \int_S^T \delta_{B_s} \, ds$.
  
Let us fix $\epsilon > 0$, $n \in N$, 
$\vec{u} \in \UU(n,r)$. Lemma \ref{comp-ps-champ}
associates to $(\epsilon,\vec{u})$ a full $\nu$-measure set
$A$ and a vector $\vec{z}=(z_1,\cdots,z_n)$ of points 
in $2r \Z \cap I_{T/\log(T)}$, such that when 
$ v \in A$, and $T$ is sufficiently large,
\begin{equation}
\label{comp-ps-v-n}
\forall i \in \{1,\cdots,n\}\, , \,\,
\nor{\frac{v}{\sqrt{\log(T)}} - \tilde{u}_i}_{\infty,z_i +I_r}
\leq \frac{\epsilon}{6}  \, ,
\end{equation}
where $\tilde{u}_i(\cdot) \triangleq u_i(\cdot -z_i)$.

Now, let us fix $(\vec{\alpha},\vec{f}) \in \DD(n,r)$ such that
$\sum_{i=1}^n \alpha_i (u_i,f_i^2)=y$. We set
$T_0=0$, $T_i =  \sum_{j=1}^{i} \alpha_j \, T$, and
$\Delta	= T/\log(T)$.
 \[
\begin{array}{l}
|Y_T-y|   
 \leq \som{i=1}{n} \va{\frac{1}{T} 
	\int_{T_{i-1}}^{T_{i-1} + \Delta} 
		\frac{v(B_s)}{\sqrt{\log(T)}} \, ds}
\\ \hspace*{2.5cm} 
+ \som{i=1}{n} \alpha_i (1-\frac{2}{\log(T)}) 
	\va{\pare{L_{T_{i-1} + \Delta}
		  ^{T_i  - \Delta }; 
		\frac{v}{\sqrt{\log(T)}}} - (f_i^2,u_i)}
\\ \hspace*{2.5cm}
  + \som{i=1}{n} \va{\frac{1}{T} 
	\int_{T_{i} -  \Delta}^{T_{i}} 
	\frac{v(B_s)}{\sqrt{\log(T)}} \, ds}
	+ 2 \frac{|y|}{\log(T)} 
\end{array}
\] 
Therefore, for $T$ sufficiently large, $ 2 \frac{|y|}{\log(T)} \leq
\frac{\epsilon}{3}$, and   
\[
\begin{array}{l}
P_0 \cro{|Y_T-y| < \epsilon} 
\\
\geq P_0 \left[ 
\forall i \in \{1,\cdots,n\}\, , \,\,
 |B_{T_{i-1}}| \leq 1; 
 \va{\frac{1}{T}\int_{T_{i-1}}
	^{T_{i-1} + \Delta} 
			\frac{v(B_s)}{\sqrt{\log(T)}}  \, ds }
			< \frac{\epsilon}{6n}; \right.
\\ \hspace*{2.5cm}
|B_{T_{i-1} + \Delta} -z_i| \leq 1;
\tau_{z_i+I_r} \circ \theta_{T_{i-1} + \Delta}
    > T_i  - \Delta	;
\\ \hspace*{2.5cm} 
 \va{\pare{L_{T_{i-1}+ \Delta}^
	     {T_i  - \Delta} ;
		 \frac{v}{\sqrt{\log(T)}}} 
		- \pare{f_i^2,  u_i}} < \frac{\epsilon}{3};
\\ \hspace*{2.5cm} \left.
\va{\frac{1}{T}\int_{T_{i} - \Delta}
	^{T_{i} } 
			\frac{v(B_s)}{\sqrt{\log(T)}}  \, ds }
			< \frac{\epsilon}{6n}
\right]
\end{array}
\]
But , on $\acc{ \tau_{z_i+I_r} \circ 
\theta_{T_{i-1} + \Delta}  > T_i - \Delta}$,
and for $v \in A$,
\[
\begin{array}{l}
\va{ \pare{ L_{T_{i-1} + \Delta}
	^{T_i  - \Delta}, \frac{v}{\sqrt{\log(T)}}} 
		- \pare{f_i^2,  u_i}}
\\ \hspace*{1cm}
\leq \nor{\frac{v}{\sqrt{\log(T)}} - \tilde{u}_i}
	_{\infty,z_i+Ir}
	+ \va{ \pare{L_{T_{i-1} + \Delta}
		^{T_i  - \Delta} ,
		\tilde{u}_i} 
		- \pare{ \tilde{f}_i^2, \tilde{u}_i}}
\\ \hspace*{1cm}
\leq \frac{\epsilon}{6}  	
+ \va{ \pare{L_{T_{i-1} + \Delta}
		^{T_i  - \Delta} ,
		\tilde{u}_i} 
		- \pare{ \tilde{f}_i^2, \tilde{u}_i}}	
\end{array}
\]
The Markov property applied recursively at times $T_{i-1}$ yields
then 
\begin{equation} 
\label{dec-temps}
P_0 \cro{|Y_T-y| < \epsilon} 
\geq  \pro{i=1}{n} U_i \, ,
\end{equation}
where
\[
U_i = \Inf_{|z|\leq 1} 
\begin{array}[t]{l}
   P_z \left[  \va{\frac{1}{T} \int_0^{\Delta} 
		\frac{v(B_s)}{\sqrt{\log(T)}} \, ds }
			< \frac{\epsilon}{6n};
	\va{B_{\Delta}-z_i}  \leq 1 ; \right.
\\ \hspace*{2cm}
	\tau_{z_i + I_r}\circ \theta_{\Delta} 
		> \alpha_i T  - \Delta;
\\ \hspace*{2cm}
	\va{ \pare{L_{\Delta}^{\alpha_i T - \Delta} ,
		\tilde{u}_i} 
		- \pare{ \tilde{f}_i^2, \tilde{u}_i}} < \frac{\epsilon}{6};
\\ \hspace*{2cm} \left.
	\va{\frac{1}{T} \int_{\alpha_i T  - \Delta}^{\alpha_iT}
		\frac{v(B_s)}{\sqrt{\log(T)}}  \, ds }
			< \frac{\epsilon}{6n};
	\va{B_{\alpha_iT}} \leq 1 \right]
\end{array} 
\]
Now, it follows from Markov property applied successively at times 
$\alpha_i T - \Delta$ and  $\Delta$,  that 
for all $i \in \{1,\cdots,n\}$,
\begin{equation}
\label{dec-i}
 U_i \geq  V_i W_i X_i \,\, ,
\end{equation} 
with 
\[\begin{array}[t]{l}
V_i = \sous{\Inf}{|z| \leq 1}
	P_z \cro{\va{\frac{1}{T} \int_0^{\Delta} 
		\frac{v(B_s)}{\sqrt{\log(T)}} \, ds }
			< \frac{\epsilon}{6n};
	\va{B_{\Delta}-z_i} \leq 1 }  \, ,
\\ 
W_i  = \hspace{-.3cm} \sous{\Inf}{z \in z_i+I_1}
	P_z \cro{ \tau_{z_i + I_r} 
		> \alpha_i T - 2 \Delta;
	\va{ \pare{L_0^{\alpha_i T - 2 \Delta} ,
		\tilde{u}_i} 
		- \pare{\tilde{f}_i^2, \tilde{u}_i}} < \frac{\epsilon}{6}} \, ,
\\
X_i  =  \hspace{-.3cm} \sous{\Inf}{z \in z_i + I_r}
	P_z \cro{ 
	\va{\frac{1}{T} \int_0^{\Delta}
		\frac{v(B_s)}{\sqrt{\log(T)}}  \, ds }
			< \frac{\epsilon}{6n};
	\va{B_{\Delta}} \leq 1} \, .
\end{array}
\]

\vspace{.5cm}
\noindent
{\bf Estimates for  $W_i$.}
By translation invariance,
\[
W_i = \Inf_{z \in I_1}
   P_z \cro{ \tau_{I_r} > \alpha_i T - \frac{2T}{\log(T)}; 
	\va{ \pare{L_0^{\alpha_i T - \frac{2T}{\log(T)}} ,
		u_i} 
		- \pare{ f_i^2, u_i}} < \frac{\epsilon}{6} } \,.
\]
It follows then from the large deviations for the occupation measure
that for all $i \in \{1,\cdots,n\}$, 
\begin{equation} 
\label{main-part}
\limiinf{T \rightarrow \infty}
\frac{1}{T} \log W_i 
\geq - \frac{\alpha_i}{2} \nor{f_i'}_2^2 \, .
\end{equation}

\vspace{.5cm}
\noindent
{\bf Estimates for  $V_i$ and $X_i$.} 
We are now going to show that 
\begin{equation}
\label{reste}
\liminf_{T \rightarrow \infty} 
\frac{1}{T} \log V_i \geq 0\, , \,\, \mbox{ and }
\liminf_{T \rightarrow \infty} 
\frac{1}{T} \log X_i \geq 0\, .
\end{equation}
Since $V_i$ and $X_i$ are treated in the same way, we give 
only the proof for $V_i$. Let $z \in I_1$. Since 
$|z_i| \leq T/\log(T)$, we have 
\[
\begin{array}{ll}
P_z \cro{
\va{B_{T/\log(T)}-z_i} \leq 1} 
& = \int_{|y+z-z_i| \leq 1} \exp \pare{- \frac{y^2}{2 T/\log(T)}}
	\, \frac{dy}{\sqrt{2 \pi T/\log(T)}} \\
& \geq  \frac{1}{\sqrt{2 \pi T/\log(T)}} \exp \pare{
	- \frac{(T/\log(T) + 1)^2 }{2T/\log(T)} } \, .
\end{array} \, .
\]
Moreover,
\[
\begin{array}{l}
P_z \cro{\va{\frac{1}{T}\int_0^{\frac{T}{\log(T)}} 
		\frac{v(B_s)}{\sqrt{\log(T)} \, ds }}
			\geq \frac{\epsilon}{6n}}
\\ \hspace{2cm}
  \leq P_z \cro{ \tau_T \leq \frac{T}{\log(T)}}
	+ \ind_{ \frac{\max_{I_T}|v|}{\sqrt{\log(T)}} 
		\geq \frac{\epsilon \log(T)}{6n}}      
\\ \hspace{2cm}
 \leq P_0 \cro{ \tau_{T-1} \leq \frac{T}{\log(T)}}
	+ \ind_{ \frac{\max_{I_T}|v|}{\sqrt{\log(T)}} 
		\geq \frac{\epsilon \log(T)}{6n}}   
\\ \hspace{2cm}
 \leq C \exp\pare{-\frac{(T-1)^2 \log(T)}{2T}} + 
		\ind_{ \frac{\max_{I_T}|v|}{\sqrt{\log(T)}} 
		\geq \frac{\epsilon \log(T)}{6n}} \, .
\end{array}
\]
Since 
\[
\begin{array}{l}
P_z\cro{\va{\frac{1}{T} \int_0^{\frac{T}{\log(T)}} 
		\frac{v(B_s)}{\sqrt{\log(T)}} \, ds }
			< \frac{\epsilon}{6n}; 
		\va{B_{\frac{T}{\log(T)}}-z_i} < 1} 
\\ \hspace*{1cm}
\geq P_z\cro{\va{B_{\frac{T}{\log(T)}}-z_i} < 1} 
- P_z \cro{\va{\frac{1}{T} \int_0^{\frac{T}{\log(T)}}
		\frac{v(B_s)}{\sqrt{\log(T)}}  \, ds }
			\geq \frac{\epsilon}{6n}}  \, ,
\end{array}
\]
we obtain
\[
\begin{array}{l}
V_i \geq \frac{1}{\sqrt{2 \pi T/\log(T)}} \exp \pare{
	- \frac{(T/\log(T) + 1)^2 }{2T/\log(T)} }
-  C \exp\pare{-\frac{(T-1)^2 \log(T)}{2T}} 
\\ \hspace*{3cm}
- \ind_{ \frac{\max_{I_T}|v|}{\sqrt{\log(T)}} 
		\geq \frac{\epsilon \log(T)}{6n}} \, .
\end{array}
\]
By \refeq{max-gaus}, the indicator is null for $T$ sufficiently large,
and we get \refeq{reste} for $V_i$.

\vspace{.5cm}
Putting together \refeq{dec-temps}, \refeq{dec-i}, \refeq{main-part},
\refeq{reste}, and 
taking the supremum over admissible $(\vec{\alpha},\vec{f})$,  we have 
proved that $\forall r>L$, $\forall \epsilon > 0$,  $\forall n \in \N$,
$\forall \vec{u} \in \UU(n,r)$,
$\nu$-a.s.,  $\forall y \in \R$,
\[
\limiinf{T \rightarrow  \infty} \frac{1}{T} \log P_0 \cro{
 \va{Y_T -y} < \epsilon} 
\geq - \hspace{-.5cm}
 \begin{array}[t]{c} 
		\Inf \\[-4pt]
		\scriptstyle{(\vec{\alpha},\vec{f} ) \in \DD(n,r)}
	\end{array} 
	\hspace{-.3cm}
\acc{I_n(\vec{\alpha},\vec{f}): 
\sum_{i=1}^{n} \alpha_i (u_i,f_i^2)=y}  \, .
\] 
This in turn implies easily \refeq{lb-u}.
\hspace*{\fill} \rule{2mm}{2mm}

\subsection{Realizing the supremum over countably many profiles.}

  We would like now to take the supremum over functions $u_1, \cdots,u_n$
. Here,
we have to be a little careful, since the ``$\nu$-a.s'' appearing
in \refeq{lb-u} depends on the functions $u_1,\cdots,u_n$. This problem
would be overcome using the separability of $\CC(\bar{I_r})$,
if the function $K^*_r$ were continuous. This is not the
case everywhere. However, assume for a moment that we could take
the supremum over admissible functions $u_i$. We would obtain that 
$\nu$-a.s., 
\[
\begin{array}{l}
\liminf_{T \rightarrow \infty} \frac{1}{T}
\log P_0(|Y_T -y| < \epsilon)
\\ 
\geq - \sous{\Inf}{\vec{u} \in \UU(n,r)}
 \sous{\Inf}{(\vec{\alpha},\vec{f}) \in \DD(n,r)}
\acc{I_n(\vec{\alpha},\vec{f}): |\sum \alpha_i(u_i,f_i^2)-y| <\epsilon}
\\

 = - \sous{\Inf}{(\vec{\alpha},\vec{f}) \in \DD(n,r)}
	\acc{ I_n(\vec{\alpha},\vec{f}):
	\exists  \vec{u} \in \UU(n,r),
	|\sum \alpha_i(u_i,f_i^2)-y| <\epsilon}
\\ 
 = - \hspace{-.5cm}
\sous{\Inf}{(\vec{\alpha},\vec{f}) \in \DD(n,r)}
\hspace{-.3cm} 
\acc{ I_n(\vec{\alpha},\vec{f}):
\hspace{-.3cm}
\begin{array}[t]{c} \Inf \\[-4pt]
	\scriptstyle{\vec{u} \in \CC(\bar{I}_r)^n}
\end{array} 
\hspace{-.5cm}
\acc{\Max_i K^*_r(u_i): |\sum \alpha_i (u_i,f_i^2)-y| < \epsilon } <1}
\end{array} 
\]

We are thus led to show that the infimum of
$\Max_i K^*_r(u_i)$ on the set $\acc{ \vec{u} \in \CC(\bar{I}_r)^n
: |\sum \alpha_i (u_i,f_i^2) -y|<\epsilon}$
can actually be reached on a countable subset of $\CC(\bar{I}_r)^n$.

\begin{Lemme}.
\begin{itemize}
\item $\forall f \in L^2(I_r)$, $K \star f^2 \in \CC(\bar{I_r})$, and
\begin{equation}
\label{k*ok}
 K^*_r(K \star f^2) = \frac{1}{2} (K \star f^2,f^2) \, .
\end{equation} 
\item
$\forall n \in \N$, $\forall (\vec{\alpha},\vec{f} ) \in \DD(n,r)$,
and $\forall y \in \R$,
\begin{equation}
\label{id-inf}
\begin{array}[t]{c} 
	\Inf \\[-4pt]
	\scriptstyle{\vec{u}\in \in \CC(\bar{I}_r)^n}
\end{array}
\hspace{-.3cm} \acc{\Max_i K^*_r(u_i): \sum \alpha_i (u_i,f_i^2) = y}
= \scriptstyle{\frac{|y|^2}
{2 \pare{\som{i=1}{n} \alpha_i \sqrt{(K \star f_i^2,f_i^2)}}^2}} \, ,
\end{equation}
with the convention $0/0=0$. Moreover, the infimum in \refeq{id-inf}
is reached for functions $(\bar{u}_1,\cdots, \bar{u}_n)$ 
defined in the following way.
Let $I_0 =\{i; \alpha_i=0 \}$. 
\begin{itemize}
\item If $\som{i=1}{n}  \alpha_i \sqrt{(K \star f_i^2,f_i^2)}=0$, take
\[
\left\{ 
\begin{array}{ll} 
	\bar{u}_i \equiv 0, &  \mbox{ for  } i \in I_0; \\ 
	\bar{u}_i \equiv \frac{y}{ \alpha_i |I_0|},  
	& \mbox{ for  } i \notin I_0 \, .  
\end{array}  
\right.  
\]
\item If $\som{i=1}{n} \alpha_i \sqrt{(K \star f_i^2,f_i^2)}>0$, take
\[
\left\{ 
\begin{array}{ll}
	\bar{u}_i \equiv 0, &  \mbox{ for  } i \in I_0; \\
	\bar{u}_i= \frac{y}
			{\sum \alpha_j \sqrt{(K \star f^2_j,f^2_j)}} 
	     \frac{K \star f^2_i}{\sqrt{(K \star f_i^2,f_i^2)}},
		& \mbox{ for } i \notin I_0.
\end{array}
\right.
\]
\end{itemize}  
\item 
Let $\D_1$ be a dense countable subset of $L^2(I_r)$, and let 
\[
\D \triangleq
\{ s (K \star g^2): g \in \D_1, s \in \{-1;1\} \}
\cup \{u\equiv q, q \in \Q\} \, .
\]
$\D$ is a countable subset of  $\CC(\bar{I_r})$, 
and $\forall n \in \N$, $\forall (\vec{\alpha}, \vec{f}) \in \DD(n,r)$,
 $\forall y \in \R$, 
$\forall  \epsilon > 0$,  
\begin{equation} 
\label{inf-denb}
\begin{array}{l}
\begin{array}[t]{c} 
	\Inf \\[-4pt]
	\scriptstyle{\vec{u}\in \CC(\bar{I}_r)^n}
	\end{array} 
\hspace{-.3cm}\acc{\Max_i K^*_r(u_i):|\sum \alpha_i (u_i,f_i^2) -y|<\epsilon}
\\
\hspace{1cm}
= \begin{array}[t]{c} 
	\Inf \\[-4pt]
	\scriptstyle{\vec{u} \in \D^n}
	\end{array}
\hspace{-.3cm} \acc{\Max_i K^*_r(u_i):|\sum \alpha_i (u_i,f_i^2) -y|<\epsilon}
\end{array}
\end{equation} 
\end{itemize}
\end{Lemme}

\noindent
{\bf Proof}.\\
{\bf Proof of \refeq{k*ok}}.
\[
\begin{array}{ll}
K^*_r(K \star f^2)
& = \Sup 
\acc{ (\mu,K\star f^2) - \frac{1}{2} (K \star \mu, \mu):
\mu \in \MM(I_r)} \\
& = \Sup 
 \acc{ \frac{1}{2} (K \star f^2,f^2) - \frac{1}{2} (K \star \mu,\mu):
\mu \in \MM(I_r)}
\, , \end{array} 
\]
by the change of variable $\mu \rightarrow \mu + f^2 \, dx$. Thus
$K^*_r(K \star f^2) = \frac{1}{2} (K \star f^2,f^2)$. 

\vspace{.5cm}
\noindent
{\bf Proof of \refeq{id-inf}}.
First of all, note that $\sum \alpha_i (\bar{u}_i, f_i^2) =y$, so that
\[ 
\begin{array}{l}
\begin{array}[t]{c} 
	\Inf \\[-4pt]
	\scriptstyle{\vec{u} \in \CC(\bar{I}_r)^n} 
\end{array}
\hspace{-.3cm} \acc{\Max_i K^*_r(u_i):\sum \alpha_i (u_i,f_i^2) = y}
\\ \hspace{1cm}
\leq \Max_i K^*_r(\bar{u}_i) = \sous{\Max}{i \notin I_0}  K^*_r(\bar{u}_i) \, .
\end{array} 
\]
If $\sum \alpha_i \sqrt{(K \star f_i^2,f_i^2)} =0$, then
$(K \star f_i^2,f_i^2)=0$
for any $i\notin I_0$. In this situation, 
$K^*_r(1) \geq \frac{(1,f_i^2)}{2 (K \star f_i^2, f_i^2)} = +\infty$, for
any $ i \notin I_0$. Thus, 
\[ \sous{\Max}{i \notin I_0}  K^*_r(\bar{u}_i) 
= \left\{ \begin{array}{ll}
 0 & \mbox{ if } y=0 \\
 + \infty & \mbox{ if } y \neq 0 
\end{array}
\right.
= \frac{y^2}{2 (\sum \alpha_i \sqrt{(K \star f_i^2,f_i^2)})^2} \, .
\]
If $\sum \alpha_i \sqrt{(K \star f_i^2,f_i^2)}> 0$, then
\[
\begin{array}{ll}
\sous{\Max}{i \notin I_0}  K^*_r(\bar{u}_i) 
& = \sous{\Max}{i \notin I_0} 
\frac{y^2}{(\sum \alpha_i \sqrt{(K \star f_i^2,f_i^2)})^2}
\frac{K^*_r(K \star f_i^2)}{(K \star f_i^2,f_i^2)}
\\
& = \frac{y^2}{2(\sum \alpha_i \sqrt{(K \star f_i^2,f_i^2)})^2} \, ,
\mbox{ by  \refeq{k*ok}.}
\end{array}
\]

\vspace{.5cm}
\noindent 
It remains now to show that
\[\begin{array}[t]{c} 
	\Inf \\[-4pt]
	\scriptstyle{\vec{u} \in \CC(\bar{I}_r)^n} 
\end{array}
\hspace{-.3cm} \acc{\Max_i K^*_r(u_i):\sum \alpha_i (u_i,f_i^2) =y}
\geq \frac{y^2}{2 (\sum \alpha_i \sqrt{(K \star f_i^2,f_i^2)})^2} \, .
\]
First, note that 
\[
\begin{array}{l}
\begin{array}[t]{c} 
	\Inf \\[-4pt]
	\scriptstyle{\vec{u} \in \CC(\bar{I}_r)^n} 
\end{array}
\hspace{-.3cm}\acc{\Max_i  K^*_r(u_i):
	\sum \alpha_i (u_i,f_i^2) =y}
\\ \hspace{1cm}
= \begin{array}[t]{c} 
	\Inf \\[-4pt]
	\scriptstyle{\vec{y} \in \R^n, \sum \alpha_i y_i =y}
\end{array}
\begin{array}[t]{c} 
	\Inf \\[-4pt]
	\scriptstyle{\vec{u}\in \CC(\bar{I}_r)^n}
\end{array}
\hspace{-.3cm} \acc{\Max_i K^*_r(u_i): \forall i, (u_i,f_i^2) =y_i }
\\ \hspace{1cm}
\geq 
\begin{array}[t]{c} 
	\Inf \\[-4pt]
	\scriptstyle{\vec{y} \in \R^n, \sum \alpha_i y_i =y}
\end{array}
\Max_i \Inf\acc{ K^*_r(u_i): u_i \in \CC(\bar{I}_r), (u_i,f_i^2) =y_i }
\end{array}
\]
For $(u_i,f_i^2)=y_i$, 
$K^*_r(u_i) \geq \frac{(u_i,f_i^2)^2}{2 (K \star f_i^2,f_i^2)} = 
\frac{y_i^2}{2 (K \star f_i^2,f_i^2)}$, so that 
\[ 
\begin{array}{l}
\begin{array}[t]{c} 
	\Inf \\[-4pt]
	\scriptstyle{\vec{u} \in \CC(\bar{I}_r)^n}
\end{array}
\hspace{-.3cm} \acc{\Max_i K^*_r(u_i):\sum \alpha_i (u_i,f_i^2) =y}
\\
\hspace{1cm}
\geq 
\begin{array}[t]{c} 
	\Inf \\[-4pt]
	\scriptstyle{\vec{y} \in \R^n}
\end{array}
\hspace{-.3cm}
\acc{
\Max_i \frac{y_i^2}{2 (K \star f_i^2,f_i^2)}: \sum \alpha_i y_i =y} \, .
\end{array}
\]
Now, for $\sum \alpha_i y_i =y$, 
\[
|y| \leq  \Max_i\cro{\frac{|y_i|}{ \sqrt{(K \star f_i^2,f_i^2)}}}
\,\,	\sum  \alpha_i \sqrt{(K \star f_i^2,f_i^2)} \, .
\]
Thus, 
\[
\begin{array}[t]{c} 
	\Inf \\[-4pt]
	\scriptstyle{\vec{u} \in \CC(\bar{I}_r)^n}
\end{array}
\hspace{-.3cm} \acc{\Max_i K^*_r(u_i):\sum \alpha_i (u_i,f_i^2) =y}
\geq \frac{y^2}{2 (\sum  \alpha_i \sqrt{(K \star f_i^2,f_i^2)})^2}
\, .
\]
This ends the proof of \refeq{id-inf}.

\vspace{.5cm}
\refeq{inf-denb} is a straightforward consequence of $\refeq{id-inf}$,
of the expression of the minimizing functions $\bar{u}_i$, and of the
continuity of $f \in L^2(I_r) \mapsto (K \star f^2,f^2)$, and
$f \in L^2(I_r) \mapsto K \star f^2 \in \CC(\bar{I}_r)$.
\hspace*{\fill} \rule{2mm}{2mm}

\vspace{1cm}

Performing now in lemma \ref{lb-u-lem} the supremum over functions
$u_i \in \D$, then over $r \in \Q$, we have thus shown
that when $K$ has compact support,  $\nu$-a.s., 
$\forall n \in \N$,
 $\forall \epsilon >0$, $\forall y \in \R$,
\begin{equation}	
\label{lb-n}
\liminf_{T \rightarrow \infty} \frac{1}{T} 
\log P_0 (|Y_T - y| < \epsilon) 
\geq - \II_n(y) \, .
\end{equation} 
where 
\begin{equation}
\label{lb-n-def}
\begin{array}{rcl}
\II_n(y) & \triangleq &
	\Inf \acc{I_n(\vec{\alpha},\vec{f}):
		(\vec{\alpha},\vec{f}) \in \DD_n(y)}
\\
\DD_n(y) & \triangleq &
\acc{ (\vec{\alpha},\vec{f}) \in \DD(n) ;
	  \frac{|y|}{\sqrt{2}} 
		<\sum \alpha_i \sqrt{(K\star f_i^2,f_i^2)}}
\, .
\end{array}
\end{equation} 

\subsection{ Identifying the rate function.}
Now, our aim is to characterize the limit $n \rightarrow \infty$
in \refeq{lb-n}. 

\begin{Lemme}.
\label{prop-In}
\begin{enumerate}	
\item $\forall n \in \N$, $\forall y \in \R$,
\begin{equation}
\label{In-dec}
\JJ(y) \leq \II_{n+1}(y) \leq \II_n(y) \, .
\end{equation}	
\item $\forall n \in \N$, $\forall \alpha \in [0,1]$, 
$\forall y_1,y_2 \in \R$,
\begin{equation}
\label{In-conv}
\II_{2n}(\alpha y_1 + (1-\alpha) y_2) 
\leq \alpha \II_n(y_1) + (1-\alpha) \II_n(y_2) \, .
\end{equation}	
\item If $\II_1^*$ denotes the Fenchel-Legendre transform 
of $\II_1$, $\II_1^{**} = \JJ$.
\item Let $\II(y) \triangleq \hspace{-6pt}
\sous{\lim}{n \rightarrow \infty} \hspace{-10pt} \searrow 
\II_n(y)$,
and $\tilde{\II}(y) \triangleq \hspace{-3pt} 
 \sous{\Sup}{\epsilon > 0} \hspace{-5pt}
\sous{\Inf}{z;|z-y| \leq
\epsilon} \hspace{-5pt} \II(z)$ the greater l.s.c. minorant of $\II$. Then 
$\tilde{\II}=\JJ$. 

\end{enumerate} 
\end{Lemme}

\noindent
{\bf Proof of 1.}
From the large deviations upper bound, we have 
$\JJ(y) \leq \II_n(y)$ for all $n$.

For any $(\vec{\alpha}, \vec{f}) \in \DD_n(y)$,
$\vec{\beta}\triangleq (\vec{\alpha},0)$ and $\vec{g}\triangleq (\vec{f},f_1)$
are such that $(\vec{\beta},\vec{g}) \in \DD_{n+1}(y)$, so that
$\II_{n+1}(y) \leq I_{n+1}(\vec{\beta},\vec{g})= I_n(\vec{\alpha},\vec{f})$.
Taking the infimum over $\DD_n(y)$ yields $\II_{n+1}(y) \leq \II_n(y)$.

\vspace{.5cm}
\noindent
{\bf Proof of 2.}
In the same way, let $\alpha \in [0,1]$ and $y_1,y_2 \in \R$ be
fixed. For any $(\vec{\beta},\vec{f}) \in \DD_n(y_1)$, and
any $(\vec{\gamma },\vec{g}) \in \DD_n(y_2)$, 
$\vec{\lambda} \triangleq (\alpha \vec{\beta}, (1-\alpha) \vec{\gamma})$
and $\vec{h} \triangleq  (\vec{f},\vec{g})$ are such that
 $(\vec{\lambda},\vec{h})
\in \DD_{2n}(\alpha y_1+(1-\alpha)y_2)$. Thus,
\[
\II_{2n}(\alpha y_1+(1-\alpha)y_2) 
\leq I_{2n}(\vec{\lambda},\vec{h}) 
= \alpha I_n(\vec{\beta},\vec{f}) + (1-\alpha) I_n(\vec{\gamma},\vec{g})
\,.
\]
Taking the infimum over elements of $\DD_n(y_1)$ and $\DD_n(y_2)$, leads
to \refeq{In-conv}.

\vspace{.5cm}
\noindent
{\bf Proof of 3.}
 Let us now compute the Le\-gen\-dre transform of $\II_1$. Fist of all,
note that
\[
\II_1(y) = \sous{\inf}{f \in H^1, \nor{f}_2=1}
	\acc{\frac{1}{2} \nor{f'}_2^2;  |y| < \sqrt{2 (K\star f^2,f^2)} } \,.
\]
Therefore,
\[
\begin{array}[t]{ll}
\II_1^*(\alpha)
& = \sous{\sup}{y \in \R} \acc{ \alpha y - \II_1(y) } \\
& = \sous{\sup}{f \in H^1, \nor{f}_2=1}
	\sous{\sup}{y \in \R} 
	\acc{ \alpha y - \frac{1}{2} \nor{f'}_2^2: 
	|y| < \sqrt{2 (K\star f^2,f^2)} } \\
& = \sous{\sup}{f \in H^1, \nor{f}_2=1} 
	\acc{ |\alpha| \sqrt{2 (K\star f^2,f^2)}  
	 - \frac{1}{2} \nor{f'}_2^2 } \\
& = \Lambda(\alpha)
\end{array}
\]
Hence, $\II_1^{**}(y)=\Lambda^*(y) =\JJ(y)$.

\vspace{.5cm}
\noindent
{\bf Proof of 4.}
Taking the limit in \refeq{In-dec}, we obtain that for
all $y \in \R$, $\JJ(y) \leq \II(y) \leq \II_1(y)$. Since
$\JJ$ is l.s.c, we also have $\JJ(y) \leq \tilde{\II}(y) \leq 
 \II(y) \leq \II_1(y)$. Since $\II_1^{**} = \JJ$, the preceding
inequality implies that $\JJ(y) = \tilde{\II}^{**}$. Now, taking
the limit in \refeq{In-conv}, we see that $\II$ is convex, and so
is $\tilde{\II}$. $\tilde{\II}$ being convex and l.s.c., 
$\tilde{\II}=\tilde{\II}^{**}=\JJ$.
\hspace*{\fill} \rule{2mm}{2mm}

\begin{Lemme}. Assume that $K$ has compact support. Then,
$\nu$-a.s., $\forall y \in \R$,
\[
\lim_{\epsilon \rightarrow 0} 
\liminf_{T \rightarrow \infty} 
\frac{1}{T} \log P_0 \cro{|Y_T-y| \leq \epsilon} \geq - \JJ(y) \, .
\]
\end{Lemme}
 
\noindent
{\bf Proof}. Taking the limit $n \rightarrow \infty$ in \refeq{lb-n}
yields that $\nu$-a.s, $\forall \epsilon >0$, $\forall y \in \R$,
\[
\liminf_{T \rightarrow \infty} 
\frac{1}{T} \log P_0 \cro{|Y_T-y| \leq \epsilon} \geq - \II(y) \, .
\]
Let $z$ be any point in $\BB(y,\epsilon)$, and let $\eta >0$ be
such that $\BB(z,\eta) \subset \BB(y,\epsilon)$. 
\[
\liminf_{T \rightarrow \infty} 
\frac{1}{T} \log P_0 \cro{|Y_T-y| \leq \epsilon}
\geq \liminf_{T \rightarrow \infty} 
\frac{1}{T} \log P_0 \cro{|Y_T-z| \leq \eta} 
\geq - \II(z) \, .
\]
Taking the supremum in $z \in \BB(y,\epsilon)$, and letting
$\epsilon$ go to 0, leads to 
\[
\lim_{\epsilon \rightarrow 0}
\liminf_{T \rightarrow \infty} 
\frac{1}{T} \log P_0 \cro{|Y_T-y| \leq \epsilon} 
\geq - \tilde{\II}(y) = - \JJ(y) 
\, .
\]
\hspace*{\fill} \rule{2mm}{2mm}

\subsection{The general case.} 
We are now going to prove
the lower bound in the general case, i.e. 
under assumption \refeq{densiteL1} for the covariance $K$.
To this end, we use the decomposition of $v = v_L+\tilde{v}_L$ 
(cf  section \ref{champ} and equation \refeq{dec-champ}). Let 
$Y=Y_L + \tilde{Y}_L$ the corresponding decomposition of $Y$.
Let $\epsilon >0$ and $L$ sufficiently large so that $\sqrt{2 \tilde{K}_L(0)}
< \frac{\epsilon}{2}$. Then, 
\[
P_0 \cro{ |Y_T - y| < \epsilon}
 \geq P_0 \cro{|Y_{L,T} - y| < \epsilon/2} 
	- P_0 \cro{|\tilde{Y}_{L,T}| \geq \epsilon/2} \,\,.
\]
But, 
\[
\begin{array}{ll}
P_0 \cro{|\tilde{Y}_{L,T}| \geq \epsilon/2} 
&  \leq P_0 \cro{ \tau_{RT} > T;|\tilde{Y}_{L,T}| \geq \epsilon/2}
+ P_0 \cro{ \tau_{RT} \leq T} 
\\
& \leq \ind_{\frac{\max_{I_{RT}} |\tilde{v}_L|}{\sqrt{\log(T)}}
	 \geq \epsilon/2} +   P_0 \cro{ \tau_{RT} \leq T} \, .
\end{array}
\]
Thus, $\forall \epsilon  >0$, and $L$ sufficiently large, $\nu$-a.s.,
\[
\liminf_{T \rightarrow \infty} 
\frac{1}{T} \log P_0 \cro{|\tilde{Y}_{L,T}| \geq \epsilon/2}
\geq - \infty.
\]
Therefore, by \refeq{lb-n}, $\forall \epsilon > 0$, $\forall L$ sufficiently
large, $\nu$.a.s, $\forall y \in \R$, $\forall n \in \N$ 
\begin{equation}
\label{lb-nL}
\liminf_{T \rightarrow \infty} \frac{1}{T} \log
P_0 \cro{ |Y_T - y| < \epsilon}
\geq - \II_n^L(y) \, ,
\end{equation} 
where 
$\begin{array}[t]{rl}
\II^L_n(y) & \triangleq 
\begin{array}[t]{c}  
	\Inf \\[-4pt] 
   	\scriptstyle{(\vec{\alpha},\vec{f}) \in \DD^L_n(y)}
\end{array} 
I_n(\vec{\alpha},\vec{f})  \, , 
\\
\DD^L_n(y) & \triangleq 
\acc{ (\vec{\alpha},\vec{f}) \in \DD(n) ;
\frac{|y|}{\sqrt{2}}  	<\sum \alpha_i \sqrt{(K_L\star f_i^2,f_i^2)}
}\, . 
\end{array} 
$

We are now going to prove that $\forall n$ and $\forall y$,
$\limsup_{L \rightarrow \infty, L \in \Q}
\II^L_n(y) \leq \II_n(y)$,     
 and we can assume that $\II_n(y) < \infty$.   
Let $\eta > 0$, and $(\vec{\alpha},\vec{f}) \in \DD_n(y)$ be
 such that $I_n(\vec{\alpha},\vec{f}) \leq {\II_n}(y)  
+ \eta$.   Since $K_L$ converges almost everywhere to $K$ when $L \rightarrow  
\infty$, $\forall i$,
$(K_L \star f_i^2,f_i^2) \sous{\rightarrow}{L \rightarrow \infty}    
(K \star f_i^2,f_i^2)$   by Lebesgue dominated convergence theorem.        
Thus, for $L$            
sufficiently large,   $(\vec{\alpha},\vec{f}) \in \DD_n^L(y)$, and       
$\II_n^L(y) \leq I_n(\vec{\alpha},\vec{f}) \leq \II_n(y) + \eta$.       
Therefore, letting first $L \rightarrow \infty$, then    
$n \rightarrow \infty$  in \refeq{lb-nL}, we obtain    
that $\nu$-a.s., $\forall \epsilon > 0$, $\forall y \in \R$,   
\[
\liminf_{T \rightarrow \infty} \frac{1}{T} \log 
P_0 \cro{ |Y_T - y| < \epsilon}   
\geq - \II(y) \, .
\]
As usual, this in turn implies the same bound with $\tilde{\II}$ in  
place of $\II$. To conclude the proof of \refeq{lby-eq},   
note that the results of  lemma \ref{prop-In} are independent  
of the support of $K$, so that we have 
$\tilde{\II}=\JJ$.

\subsection{ Lower bound for $X_{2,T}$}
As for the upper bound, the lower bound for $Y_T$ yields 
straightforwardly the same lower bound for $X_{2,T}$, since
\[
P_0 \cro{ \va{\frac{X_{2,T}}{T \sqrt{\log(T)}} - y} < \epsilon } 
\geq P_0 \cro{|Y_T -y | < \epsilon/2} 
-P_0 \cro{ \va{\frac{W_{2,T}}{T \sqrt{\log(T)}}} \geq \epsilon/2} \, ,
\]
and $\lim_{T \rightarrow \infty} \frac{1}{T} \log 
P_0 \cro{ \va{\frac{W_{2,T}}{T \sqrt{\log(T)}}} \geq \epsilon/2} 
= - \infty$.

\vspace{1cm}
\section{Properties of the rate function.}
\label{rate}

The aim of this section is to prove proposition \ref{prop-J} linking
the behavior
of $K$ at infinity, with the behavior of  $\JJ$ near
the origin.  Note
that since the functions normalizing $\JJ$ are convex and continuous,
and since $\II_1^{**}=\JJ$, 
to prove  proposition \ref{prop-J}, it is
enough to prove the corresponding assertions for $\II_1$.

Using the isometry 
of $L^2$: $f \mapsto f_{\lambda} =\sqrt{\lambda} f(\lambda \cdot)$,
note that 
\begin{equation} 
\label{isometrie}
\begin{array}{ll}
\II_1(y) 
& = \hspace{-5mm} \sous{\Inf}{f \in H^1, \nor{f}_2=1} \hspace{-3mm}
	\acc{  \frac{\lambda^2}{2} \nor{f'}_2^2;  \,
	\int \hspace{-.2cm} \int 
	K (\frac{x-z}{\lambda}) f^2(x) f^2(z) \, dx \, dz 
	> \frac{y^2}{2} } , \, \forall \lambda >0 \\
& = \hspace{-5mm}
	\sous{\Inf}{f \in H^1, \nor{f}_2=1}  \hspace{-2mm}
    \sous{\Inf}{\lambda > 0}  \hspace{-2mm}
	\acc{  \frac{\lambda^2}{2} \nor{f'}_2^2;  \,
	 \int \hspace{-.2cm} \int 
	K \pare{\frac{x-z}{\lambda}} f^2(x) f^2(z) \, dx \, dz 
	> \frac{y^2}{2} }
\end{array}
\end{equation}

\noindent  
{\bf Case $\limsup_{|x| \rightarrow \infty} |K(x)| |x|^{\beta} < \infty$ 
for some $\beta \in ]0,1[$.} \\
Since $K$ is bounded, there exists a constant $C$ such that
$K(x) \leq C |x|^{-\beta}$. It follows then from  \refeq{isometrie}
that
\[ 
\II_1(y)
 \geq \sous{\Inf}{f \in H^1(\R), \nor{f}_2=1} 
\sous{\Inf}{\lambda >0}
\acc{\frac{\lambda^2}{2} \nor{f'}_2^2:
C \lambda^{\beta} (I_{\beta}(f^2),f^2) > \frac{y^2}{2}}
 \, ,
\]
where $I_{\beta}$ is the Riesz operator defined 
by
$I_{\beta}(f)(x) \triangleq \int_{\R} \frac{f(y)}{|x-y|^{\beta}} \, dy$.

Taking the infimum in $\lambda$ 
leads to 
\begin{equation}
\label{ub-I1}
\frac{\II_1(y)}{|y|^{\frac{4}{\beta}}} 
 \geq C \Inf 
\acc{\frac{\nor{f'}_2^2}{ (I_{\beta}(f^2),f^2)^{\frac{2}{\beta}}}:
f \in H^1(\R), \nor{f}_2=1} \, .
\end{equation} 

 Now, for $p \in ]1, \frac{1}{1-\beta}[$, $I_{\beta}$ is continuous from
$L^p(\R)$ to $L^r(\R)$, for $1/r=1/p -(1-\beta)$ (see for instance
theorem 1 pp 119 in \cite{stein}).
Therefore, for any $f \in H^1(\R)$ such that $\nor{f}_2=1$,
and for any $p \in ]1, \frac{1}{1-\beta}[$,
$$\begin{array}{lll}
(I_{\beta}(f^2),f^2) & \leq \nor{f^2}_{r'} \nor{I_{\beta}(f^2)}_r 
	& \mbox{ where } \Frac{1}{r} + \Frac{1}{r'} = 1 \, , \\ 
& \leq C \nor{f^2}_{r'} \nor{f^2}_p  
	& \mbox{ by continuity of } I_{\beta} \, ,\\
& \leq C \nor{f}_{\infty}^{2(1-\frac{1}{r'}+1-\frac{1}{p})} 
	& \mbox{ since } \int f^2 =1 \, .
\end{array}
$$
Note that  by Sobolev embedding theorem, any function $f \in H^1(\R)$
belongs to  $L^{\infty}(\R)$, and 
$\nor{f}_{\infty} \leq C \nor{f}_2^{1/2} \nor{f'}^{1/2}_2$
for some cons\-tant $C \in ]0,\infty[$. Thus,
for any $f \in H^1(\R)$ such that $\nor{f}_2=1$, 
\[
(I_{\beta}(f^2),f^2) \leq C \nor{f'}_2^{\beta} \,.
\]
Therefore, the infimum in \refeq{ub-I1} 
is strictly positive, and it is clearly finite.

\vspace{.5cm}
Let us now turn to the converse inequality, and let us assume 
that $K \geq 0$ and $l= \liminf_{|x| \rightarrow \infty} K(x) |x|^{\beta} > 0$.
The change of variable $\lambda = \gamma |y|^{\frac{2}{\beta}}$
in \refeq{isometrie} leads to 
\[
\frac{\II_1(y)}{|y|^{\frac{4}{\beta}}}
=    \hspace{-.5cm}
\sous{\Inf}{f \in H^1(\R), \nor{f}_2=1} 
\hspace{-.3cm}     \sous{\Inf}{\gamma >0}
	\acc{ \frac{\gamma^2}{2} \nor{f'}_2^2:
 \frac{1}{|y|^2} 
	\int \hspace{-.3cm} \int 
	K (\frac{x-z}{\gamma |y|^{2/\beta}}) 
	f^2(x) f^2(z) dx dz > \frac{1}{2}}   \, .
\]
Let $f \in H^1(\R)$, $\nor{f}_2=1$, and $\gamma > 0$ be such that
$l \gamma^{\beta} (I_{\beta}(f^2),f^2) > \frac{1}{2}$. By Fatou lemma,
\[
l \gamma^{\beta} (I_{\beta}(f^2),f^2) \leq 
\limiinf{|y| \rightarrow 0}  
 \frac{1}{|y|^2} 
	\int \int K \pare{\frac{x-z}{\gamma |y|^{\frac{2}{\beta}}}} 
	f^2(x) f^2(z) \, dx \, dz \, ,
\]
and thus, for any $(f,\gamma)$ with $\nor{f}_2=1$, and 
$l \gamma^{\beta} (I_{\beta}(f^2),f^2)  > \frac{1}{2}$, 
\[
\limsup_{|y| \rightarrow 0}  
\frac{\II_1(y)}{|y|^{\frac{4}{\beta}}} 
\leq \frac{\gamma^2}{2} \nor{f'}_2^2 \,
\]
Therefore,
\[
\limsup_{|y| \rightarrow 0} 
\frac{\II_1(y)}{|y|^{\frac{4}{\beta}}}
\leq C \Inf 
\acc{  \frac{\nor{f'}_2^2}{ (I_{\beta}(f^2),f^2)^{\frac{2}{\beta}}}:
f \in H^1(\R), \nor{f}_2=1}
< \infty \, .
\]

Note that using Lebesgue dominated convergence theorem in place of 
Fatou lemma, the same result holds, as soon
as $\lim_{|x| \rightarrow \infty}  K(x) |x|^{\beta}> 0$.
This concludes the proof of point 1. of proposition \ref{prop-J}.

\vspace{.5cm}
\noindent
{\bf Case $\limsup_{|x| \rightarrow \infty} |K(x)| |x|^{\beta} < \infty$
for some $\beta > 1$, $\int K(x) \, dx \neq 0$.} \\
In this situation,
\begin{equation}
\label{K-integ}
K \in L^1(\R)\, , \mbox{ and } \forall \delta  \in ]0, \beta-1[\, , \,\,\, 
\int K^2(x) |x|^{1+2\delta} \, dx < \infty \, .
\end{equation}
 Therefore,  
by dominated convergence, 
\begin{equation}
\label{lim-b>1}
 \int\int \frac{1}{\lambda} K (\frac{x-y}{\lambda}) 
f^2(x) f^2(y) \sous{\longrightarrow}{\lambda \rightarrow 0}   
\pare{\int K(x) dx}  \, \nor{f}_4^4 \,  ,
\end{equation}
Note that this implies that $\bar{K} \triangleq \int K(x) \, dx \geq 0$,
 and thus
$\bar{K} > 0$.

Moreover, we know from standard results in functional analysis 
(see for instance \cite{stein}) that 

\begin{Lemme}
\label{Sobolev}
If $f \in H^1(\R)$, then $\forall  p \in [2, +\infty]$, $f \in L^p(\R)$,
and 
\[ \nor{f}_p \leq C \nor{f}_2^{\frac{1}{2} +\frac{1}{p}}
\nor{f'}_2^{\frac{1}{2} -\frac{1}{p}} \, , 
\]
 where $C$ is a constant 
depending only on $p$. Moreover, $\forall  \delta  \in ]0,1[$,
there exists a constant $C$ such that 
\[
\pare{ \int_{\R}
\frac{\nor{f(\cdot +t) - f(\cdot)}_2^2}{|t|^{1+2\delta}} \, dt}^{\frac{1}{2}}  
\leq C \pare{\nor{f}_2+ \nor{f'}_2} \, .
\]
\end{Lemme}

\noindent 
Thus, for all $f \in H^1(\R)$, $\nor{f}_2=1$, 
 $\forall \delta  \in ]0,1[\cap ]0,\beta-1[$,  $\forall \lambda >0$,   
\begin{equation}  
\label{erreur}       
\begin{array}{l}        
\va{ \int \int   
	\frac{1}{\lambda} K \pare{\frac{x-y}{\lambda}} f^2(x) f^2(y) 
		\, d x \, d y
- (\int K(x) \, dx)  \, \int f^4(x) \, d x} \\
\leq \int dx f^2(x) \int dz |K(z)| \va{ f^2(x + \lambda z) - f^2(x)} \\
\leq 2 \nor{f}_{\infty} \int dz | K(z)| 
	\int dx f^2(x) \va{f(x+\lambda z) - f(x)} \\
\leq  2 \nor{f}_{\infty} \nor{f}_4^2  
	\int |K(z)|    
	\nor{f(\cdot + \lambda  z) - f(\cdot)}_2 \, d z \\  
\leq    2 \lambda^{\delta} \nor{f}_{\infty} \nor{f}_4^2  
	\pare{\int  K^2(z) |z|^{1+2 \delta}  
		\, d z}^{1/2}  
	\pare{  \displaystyle{\int} 
		\frac{\nor{f(\cdot - z) - f(\cdot)}_2^2}
		{|z|^{1+2 \delta}}  \, d z}^{1/2} \\ 
\leq C  \lambda^{\delta} \nor{f'}_2 (1 + \nor{f'}_2) \, \mbox{ by lemma
	\ref{Sobolev} }    \, . 
\end{array}   
\end{equation} 

\noindent
Now, the change of variable $\lambda = \gamma |y|^2$ in \refeq{isometrie}
gives
\[
\frac{\II_1(y)}{y^4}
=  \hspace{-.5cm} \sous{\Inf}{f \in H^1(\R), \nor{f}_2=1} 
   \hspace{-.3cm}  \sous{\Inf}{\gamma >0}
	\acc{ \frac{\gamma^2}{2} \nor{f'}_2^2: 
	 \frac{1}{y^2} 
	\int \hspace{-.3cm} \int K (\frac{x-z}{\gamma y^2}) 
	f^2(x) f^2(z) dx dz > \frac{1}{2}}   \, .
\]
Let us fix $(f,\gamma) \in H^1(\R) \times ]0,\infty[$ such that 
$\nor{f}_2=1$, and $\gamma \bar{K} \nor{f}_4^4 > 1/2$. \refeq{lim-b>1}
implies that 
\[
\limsup_{y \rightarrow 0} \frac{\II_1(y)}{y^4}
\leq \frac{\gamma^2}{2} \nor{f'}^2_2 \, .
\]
Taking the infimum in $\gamma$ first, then in $f$, we obtain
\[
\limsup_{y \rightarrow 0} \frac{\II_1(y)}{y^4}
\leq \frac{I}{8 \bar{K}^2} \, ;
\]
where $I \triangleq \Inf \acc{
	\frac{\nor{f'}_2^2}{\nor{f}_4^8}: f \in H^1(\R), \nor{f}_2=1}
  \in ]0,+\infty[$, by lemma
\ref{Sobolev}. 

For the opposite direction, we begin by rewriting the 
first equality in \refeq{isometrie} with $\lambda=y^2$:
\[
\frac{\II_1(y)}{y^4}
= \sous{\Inf}{f \in H^1(\R), \nor{f_2}=1}
 \acc{ \frac{1}{2} \nor{f'}_2^2: 
\frac{1}{y^2} \int \hspace{-4pt} \int 
	K(\frac{x-z}{y^2}) f^2(x) f^2(z) \, dx \, dz 
> \frac{1}{2}} \, .
\]
Let $\eta > 0$ and for each $y$, let $f_y$ satisfying the 
above constraints and $\frac{1}{2} \nor{f_y'}_2^2
	\leq \frac{\II_1(y)}{y^4} + \eta$. 
Since $\limisup{y \rightarrow 0}\frac{\II_1(y)}{y^4} < \infty$,
we also have $\limisup {y \rightarrow 0} \nor{f_y'}_2 < \infty$. 
Moreover, by \refeq{erreur}, 
\[
\va{\bar{K} \nor{f_y}_4^4 - 
\frac{1}{y^2} \int \int K(\frac{x-z}{y^2}) f_y^2(x) f_y^2(z) \, dx \, dz }
\leq C |y|^{2 \delta} \nor{f'_y}_2 (1+ \nor{f'_y}_2) \, ,
\]
so that 
$\limi{y\rightarrow 0} \va{\bar{K} \nor{f_y}_4^4 - 
\frac{1}{y^2} \int \int K(\frac{x-z}{y^2}) f_y^2(x) f_y^2(z) \, dx \, dz } =0$.
Thus, 
\[ \liminf_{y \rightarrow 0} \bar{K} \nor{f_y}_4^4 \geq \frac{1}{2} \,  .
\]
Now, by definition of $I$, $\nor{f'_y}_2^2 \geq I \nor{f_y}_4^8$.
Thus, $\liminf_{y \rightarrow 0} \frac{1}{2}  \nor{f'_y}_2^2
\geq \frac{I}{8 \bar{K}^2}$. This ends the proof of point 4. of
proposition \ref{ldp-theo}.
\vspace{0,5cm}

\noindent{\bf Acknowledgements}. We would like to thank
Francis Comets for having given us the right scaling.


\vspace{1cm}
 \begin{thebibliography}{99}

\bibitem{avellaneda-majda90} M. Avellaneda, A. Majda.
{\it Mathematical models with exact renormalization for turbulent
transport.}
Commun. Math. Phys. {\bf 131} (1990), pp 381-429.

\bibitem{avellaneda-majda92} M. Avellaneda, A. Majda.
{\it Mathematical models with exact renormalization for turbulent
transport II.}
Commun. Math. Phys. {\bf 146} (1992), pp 139-204.

\bibitem{castell-pradeilles} F. Castell, F. Pradeilles.
{\it Annealed large deviations for diffusions in a random
Gaussian shear flow drift.}
Stoc. Proc. and Appl. {\bf 94} (2001), pp 171-197.


\bibitem{carmona} R. Carmona,
{\it Transport properties of Gaussian velocity fields.}
in Real and Stochastic Analysis. Probab. Stochastics Series.
CRC, Boca Raton, FL (1997), pp 9-63.

\bibitem{carmona-xu} R. Carmona, L. Xu,
{\it Homogenization for time dependent 2-D incompressible Gaussian flows.}
Ann. Appl. Probab. {\bf 7} (1997), no 1, pp 265-279.







 

\bibitem{gartner-konig} J. G\"artner, W. K\"onig.
{\it Moment asymptotics for the continuous parabolic Anderson model.}
Ann. Appl. Probab. {\bf 10} (2000), no 1, 192--217.

\bibitem{gartner-konig-molchanov} J. G\"artner, W. K\"onig, S. A. Molchanov.
{\it Almost sure asymptotics for the continuous parabolic Anderson
model.}
Probab. Theory Related  Fields {\bf 118} (2000), no 4, 547--573.

 
\bibitem{landim-olla-yau} C. Landim, S. Olla, H.T. Yau.
{\it Convection-diffusion equation with space-time ergodic random flow.}
Probab. Theory Related Fields {\bf 112} (1998), no 2, 203--220.
 


\bibitem{olla} S. Olla.
{\it Homogenization of diffusion processes in random fields.}
Cours de l'Ecole Polytechnique (1994).

\bibitem{stein} E. M. Stein.
{\it Singular integral and differentiability properties of functions.}
Princeton Mathematical series. No 30. Princeton University Press,
Princeton, N.J. 1970. 

\bibitem{sznitman} A. S. Sznitman.
{\it Brownian motion, obstacles and random media.}
Sprin\-ger Monographs in Mathematics. Springer-Verlag, Berlin, 1998. 




\end{thebibliography}
\end{document}